\documentclass[12pt,fleqn]{article}

\usepackage{enumitem} 
\usepackage{a4}
\usepackage{amsthm}
\usepackage{ifthen}
\usepackage{amsmath}
\usepackage{amssymb}
\usepackage{mathtools}
\usepackage{url}
\usepackage{hyphenat} 

\swapnumbers

\theoremstyle{plain}
\newtheorem{theorem}{Theorem}[section]
\newtheorem{lemma}[theorem]{Lemma}
\newtheorem{proposition}[theorem]{Proposition}

\theoremstyle{definition}
\newtheorem{definition}[theorem]{Definition}

\newtheorem{remark}[theorem]{Remark}
\newtheorem{notation}[theorem]{Notation}

\numberwithin{equation}{theorem} 

\newcommand{\New}[1]{\emph{#1}}

\DeclareMathOperator{\Pol}{Pol}
\newcommand{\Pola}[1][]{\Pol^{(#1)}} 
\DeclareMathOperator{\Inv}{Inv}
\DeclareMathOperator{\Op}{Op}
\newcommand{\Opa}[1][]{\Op^{(#1)}}
\DeclareMathOperator{\Rel}{Rel}

\DeclareMathOperator{\pmOp}{\prescript{\pm}{}Op} 
\newcommand{\pmOpa}[1][]{\pmOp^{(#1)}}
\DeclareMathOperator{\pmRel}{\prescript{\pm}{}Rel}

\DeclareMathOperator{\SOp}{\prescript{S}{}Op}

\newcommand{\SOpa}[1][]{\SOp^{(#1)}}

\newcommand{\SJ}{\prescript{S}{}\!\!J} 
\newcommand{\pmD}{\prescript{\pm}{}\!D} 

\newcommand{\SSg}[2][]{\ensuremath{\prescript{S}{}\!\langle #2 \rangle_{#1}}}

\newcommand{\SL}{\prescript{S}{}\!\cL}

\DeclareMathOperator{\SPol}{\prescript{S}{}Pol}

\DeclareMathOperator{\SInv}{\prescript{S}{}Inv}

\newcommand{\Spreserves}{\mathrel{\stackrel{S}{\triangleright}}}
\DeclareMathOperator{\sgn}{sgn} 
\DeclareMathOperator{\Sgn}{Sgn}

\DeclareMathOperator{\ar}{ar} 

\newcommand{\pF}[1][F]{#1^{\boxdot}} 

\DeclareMathOperator{\pr}{pr} 
 
\newcommand{\pair}[2]{\left[\begin{smallmatrix}#1\\#2\end{smallmatrix}\right]}
\newcommand{\triple}[3]{\left[\begin{smallmatrix}#1\\#2\\#3\end{smallmatrix}\right]}
\newcommand{\preserves}{\mathrel{\triangleright}}

\newcommand{\id}{{\sf id}}

\newcommand{\cL}{\mathcal{L}}
\newcommand{\cM}{\mathcal{M}}

\newcommand{\blambda}{\boldsymbol{\lambda}}

\newcommand{\bb}{\mathbf{b}}
\newcommand{\bz}{\mathbf{z}}
\newcommand{\bx}{\mathbf{x}}
\newcommand{\by}{\mathbf{y}}

\newcommand{\N}{\mathbb{N}} 

\newcommand{\clUinf}{\mathsf{U}^\infty}
\newcommand{\clLambda}{\mathsf{\Lambda}}

\newcommand{\Sg}[2][]{\ensuremath{\langle #2 \rangle_{#1}}}
\newcommand{\Rg}[2][]{\ensuremath{[#2]_{#1}}}

\let\rho=\varrho
\let\epsilon=\varepsilon
\let\phi=\varphi
\let\kappa=\varkappa

\newcommand{\pmRela}[1][]{\pmRel^{(#1)}}

\newcommand{\pmJ}{\prescript{\pm}{}\!\!J}

\newcommand{\pmSg}[2][]{\ensuremath{\prescript{\pm}{}\!\langle #2 \rangle_{#1}}}
\newcommand{\pmRg}[2][]{\ensuremath{\prescript{\pm}{}\![#2]_{#1}}}
\newcommand{\pmL}{\prescript{\pm}{}\!\cL}
\newcommand{\pmLsd}{\prescript{\pm}{}\!\cL^{\mathrm{sd}}}

\DeclareMathOperator{\pmPol}{\prescript{\pm}{}Pol}
\newcommand{\pmPola}[1][]{\pmPol^{(#1)}} 
\DeclareMathOperator{\pmInv}{\prescript{\pm}{}Inv}
\newcommand{\pmInva}[1][]{\pmInv^{(#1)}} 
\newcommand{\pmpreserves}{\mathrel{\stackrel{\pm}{\triangleright}}}

\newcommand{\refI}[1]{{\rm #1{${}^{\text{\tiny(I)}}$}}}

\let\OLDrestriction\restriction  
\renewcommand{\restriction}{\mathclose\OLDrestriction}

 \DeclareMathOperator{\essar}{essar} 

\author{Peter Jipsen\\ Faculty of Mathematics
  \\ Chapman University, Orange, CA (USA)
\and Erkko Lehtonen\\
Institute for Algebra \\
Johannes Kepler University Linz 
(Austria)
\and Reinhard P\"oschel%
 \\Institute of Algebra\\Faculty of Mathematics\\TU Dresden (Germany)}
\date{Nov.\ 17, 2026}

\title{$S$\hyp{}preclones and the Galois connection $\SPol$--$\SInv$, Part
  II: Boolean $\pm$\hyp{}preclones}

\begin{document}

\maketitle

\begin{abstract}
  We continue the investigation of $S$\hyp{}preclones initiated in Part~I,
  specializing to the case of Boolean $\pm$\hyp{}preclones where
  $S = \{+,-\}$ is the two\hyp{}element group and the base set is
  $A = \{0,1\}$. The main results include the complete classification
  of maximal and minimal Boolean $\pm$\hyp{}preclones in the lattice
  $\pmL_2$ of all Boolean $\pm$\hyp{}preclones. We prove that there are
  exactly nine maximal Boolean $\pm$\hyp{}preclones.

 Hence we
  obtain a completeness theorem for $\pm$\hyp{}preclones: a set of
  $\pm$\hyp{}operations generates all $\pm$\hyp{}operations if and only if it
  contains an operation that fails to preserve each of the nine
  defining $\pm$\hyp{}relations of the maximal $\pm$\hyp{}preclones. Next we
  classify the twenty minimal Boolean $\pm$\hyp{}preclones (atoms of the
  lattice), each described via a generating $\pm$\hyp{}operation.

  We also investigate self\hyp{}dual $\pm$\hyp{}preclones, characterizing the
  eight maximal and nine minimal self\hyp{}dual $\pm$\hyp{}preclones and proving
  a corresponding completeness theorem. Additionally, we include the
  result, due to Vojt\v{e}ch David, that the lattice $\pmL_2$ is
  uncountable, in contrast with the countability of the classical Post
  lattice of Boolean clones. 
This is established by embedding an
  uncountable poset of clonoids into
  $\pmL_2$.

Finally, we study preimage classes of classical Boolean
  clones within $\pmL_2$, providing complete descriptions for the
  unary case and for several non\hyp{}unary minimal clones. We also prove a
  Baker--Pixley type theorem for $\pm$\hyp{}preclones containing the
  majority operation, showing that such $\pm$\hyp{}preclones are determined
  by their binary $\pm$\hyp{}relations.
\end{abstract}

\setcounter{section}{6}
\setcounter{figure}{1}


\section*{Introduction}

An $S$\hyp{}preclone, introduced in the preceding paper, Part~I
\cite{JipLP2024}, is a set of $S$\hyp{}operations $f: A^n \to A$ where each
argument is assigned a signum $s \in S$, and $S$ is a monoid
determining the behavior of these signa under
composition. $S$\hyp{}preclones generalize classical clones and are a
subclass of preclones (operads).

In Part~I, the definitions of $S$\hyp{}operations, $S$\hyp{}preclones,
$S$\hyp{}relations, and $S$\hyp{}relational clones were established. It was
shown that these concepts are paired via the Galois connection
$\SPol$--$\SInv$ induced by the $S$\hyp{}preservation property
($f \Spreserves \varrho$). The central result of Part~I was that the
$S$\hyp{}preclones and $S$\hyp{}relational clones are exactly the Galois closed
sets when $A$ is finite (Theorems \refI{4.1} and \refI{4.6})\footnote{To maintain a unified numbering scheme across Parts I and II, the numbering of sections and all numbered mathematical statements continues from where Part~I ends. References to material from Part~I are additionally highlighted with a superscript (I).}. Part~I
also proved basic facts about the structure of the lattice of all
$S$\hyp{}preclones ($\SL_A$), demonstrating that it is atomic and coatomic,
with finitely many atoms and coatoms (Proposition \refI{5.1},
Proposition \refI{5.2}).

The specific case where $S = \{+,-\}$ is the two\hyp{}element group
equipped with multiplication was identified in Part~I as the
motivating example. In this case, $S$\hyp{}operations are referred to as
$\pm$\hyp{}operations, and the resulting structures are called
$\pm$\hyp{}preclones.

In this second part of our investigation, we specialize to the
$\pm$\hyp{}preclone case and restrict the base set to the Boolean set
$A = \{0, 1\}$. The lattice of Boolean $\pm$\hyp{}preclones is denoted by
$\pmL_2$.

A primary question concerning $\pmL_2$ is its cardinality: is the size
of $\pmL_2$ countable or uncountable (see Section \ref{sec:latticeL2})? This question
is significant because there are countably many classical Boolean
clones in Post's lattice, while in general there are uncountably many
Boolean preclones (Theorem \ref{thm:uncountable}).

This paper provides a detailed structural analysis of the lattice
$\pmL_2$. We briefly review preliminaries in Section \ref{sec:prel} and the Galois
connection $\SPol$--$\SInv$ in the setting of $\pm$\hyp{}preclones. We also
introduce the notions of $S$\hyp{}algebras and $S$\hyp{}terms required for
proofs in subsequent sections. Further details and examples about the
fundamental definitions of $S$\hyp{}preclones and $S$\hyp{}relational clones are
provided in Part~I\@.

In Section~\ref{sec:Max} we classify the maximal Boolean $\pm$\hyp{}preclones. Drawing
on results from Part~I (Proposition \refI{5.3}) characterizing maximal
$S$\hyp{}preclones, we find that there are exactly nine
maximal Boolean $\pm$\hyp{}preclones (Theorem \ref{Max-theorem}). This list enables the
formulation and proof of a completeness theorem for $\pm$\hyp{}preclones
(Theorem \ref{thm:complete1}).

We classify the minimal Boolean $\pm$\hyp{}preclones in Section \ref{sec:Min}. These
are the atoms of the lattice $\pmL_2$. We demonstrate that there are
twenty minimal Boolean $\pm$\hyp{}preclones, categorized by their minimal
signa: $(+,+,+, -)$, $(+,+,-,-)$, $(+,+,-)$, and $(+,-)$ (Theorem \ref{min-theorem}).
   
In Section \ref{sec:selfdual} we investigate self\hyp{}dual $\pm$\hyp{}preclones. We characterize
the eight maximal self\hyp{}dual $\pm$\hyp{}preclones and the nine minimal
self\hyp{}dual $\pm$\hyp{}preclones (Proposition \ref{prop:selfdual2}). Since the lattice of
self\hyp{}dual $\pm$\hyp{}preclones ($\pmLsd_{2}$) is atomic and coatomic
(Proposition \ref{prop:SDlattice}), we deduce a completeness theorem for self\hyp{}dual
$\pm$\hyp{}preclones (Theorem \ref{thm:complete2}).

We discuss further aspects of the lattice $\pmL_2$ in the final Section~\ref{sec:latticeL2},
including its cardinality (Theorem \ref{thm:uncountable} due to Vojt\v{e}ch David),
the concept of preimage classes ($\pF{}$) for
$F \in \mathcal{L}_2$,
and Baker--Pixley theorem for $\pm$\hyp{}preclones (Proposition \ref{prop:BakerPixley}).


\section{Preliminaries}
\label{sec:prel}

In this section, we first briefly revisit the basic concepts and
notations pertaining to $S$\hyp{}preclones, presenting these definitions in
a specialized form for so\hyp{}called $\pm$\hyp{}preclones. For further details
and examples, we refer the reader to Part~I\@.

We let $S$ be the two\hyp{}element group. More specifically, we choose
$S = \{+1, -1\}$ equipped with multiplication, for short
$S = \{\mathord{+}, \mathord{-}\}$ (see our ``motivating example'',
Example~\refI{2.5}). Thus, $+$ is the unit element of $S$. In many
notations and terms, we use the symbol $\pm$ to refer to this choice
of $S$.

The set of all $n$\hyp{}ary operations on a set $A$ is denoted by
$\Opa[n](A)$, and $\Op(A) := \bigcup_{n \in \N_{+}} \Opa[n](A)$. A
\emph{$\pm$\hyp{}operation} is an operation $f \colon A^n \to A$ together
with a \emph{signum}
$\sgn(f) = (s_1, \dots, s_n) \in \{\mathord{+}, \mathord{-}\}^n$ that
gives labels $+$ or $-$ (\emph{signs}) to the arguments of $f$. The
set of all $n$\hyp{}ary $\pm$\hyp{}operations on $A$ is denoted by
$\pmOpa[n](A)$, and $\pmOp(A) := \bigcup_{n \in \N_{+}}
\pmOpa[n](A)$. The \emph{underlying function} of a $\pm$\hyp{}operation $f$
is the operation $\mathring{f}$ that is obtained by simply ignoring
the signum. (See Definition~\refI{2.2}.) \emph{Trivial
  $\pm$\hyp{}operations} or (\emph{trivial}) \emph{projections} are
$\pm$\hyp{}operations whose underlying operation is a projection and the
sign of the essential argument is $+$.

A \emph{$\pm$\hyp{}preclone} is a set $F \subseteq \pmOp(A)$ of
$\pm$\hyp{}operations that contains the identity operation $\id_A$ with
signum $(\mathord{+})$ and is closed under the operations $\zeta$
(cyclic shift), $\tau$ (transposition of the first two arguments),
$\nabla^s$ (adding a fictitious argument with signum
$s \in \{ \mathord{+}, \mathord{-} \}$ at the first place), $\Delta$
(identification of the first two arguments if they have the same
sign), and $\circ$ (composition). The $\pm$\hyp{}preclone generated by
$F \subseteq \pmOp(A)$ is denoted by $\pmSg{F}$ or $\pmSg[A]{F}$. The
lattice of $\pm$\hyp{}preclones on $A$ is denoted by $\pmL_{A}$. (See
Definition~\refI{2.3}.)

An $m$\hyp{}ary \emph{$\pm$\hyp{}relation} on $A$ is a pair
$\rho = (\rho_{+}, \rho_{-})$ of $m$\hyp{}ary relations
$\rho_{+}, \rho_{-} \subseteq A^m$. The set of all $m$\hyp{}ary
$\pm$\hyp{}relations on $A$ is denoted by $\pmRela[m](A)$, and
$\pmRel(A) := \bigcup_{m \in \N_{+}} \pmRela[m](A)$. (See
Definition~\refI{2.8}.) An $m$\hyp{}ary $\pm$\hyp{}relation
$\rho = (\rho_{+}, \rho_{-})$ is a \emph{$\pm$\hyp{}diagonal} if
$\rho_{+} = \rho_{-} = \delta^m_\top := \emptyset$ or
$\rho_{+} = \rho_{-} = \delta^m_\epsilon$ for some equivalence
relation $\epsilon$ on $\{1, \dots, m\}$, where $\delta^m_\epsilon$ is
the \emph{diagonal relation}
\[
  \delta^m_\epsilon := \{ \, (a_1, \dots, a_m) \in A^m \mid (i,j) \in
  \epsilon \implies a_i = a_j \, \}.
\]
(See Definition~\refI{2.9} and Remark~\refI{2.10}.)

It is customary to represent an $m$\hyp{}ary relation as a matrix with $m$
rows whose columns are the tuples belonging to the relation. In this
way, we may refer to the elements at a specific position within these
tuples collectively as a ``row'' of the relation (cf.\ also Definition \refI{2.7}). We are going to
introduce several operations on relations, and many of them can be
described intuitively in terms of simple row operations on matrices.

\begin{definition}[$\pm$\hyp{}relational clones (cf.\ Definition
  \refI{2.11}, Remark \refI{2.14(A)})]\label{rel-rules}
  A set of $\pm$\hyp{}relations is called a \New{$\pm$\hyp{}relational clone} if
  it contains $\delta^{\pm}$ and is closed under the operations
  $\zeta$, $\tau$, $\pr$, $\times$, and $\cap$ (denoted by $\land$ in Definition \refI{2.11}), as well as under the
  operation $\mu_{-}$, which we will refer to as \New{index
    translation by `$\mathord{-}$\!'} and which simply exchanges positive and
  negative parts of $\rho$. These operations are defined as follows. For
  $\rho = (\rho_{+}, \rho_{-})$ and $\rho' = (\rho'_{+}, \rho'_{-})$
  in $\pmRel(A)$, we let
  \begin{enumerate}[label=\textup{(\arabic*)}]
  \item\label{rel-R1} $\delta^{\pm} := (\Delta_{A}, \Delta_{A})$,

  \item\label{rel-R2}
    $\zeta(\rho) = \zeta ((\rho_{+}, \rho_{-})) := (\zeta \rho_{+},
    \zeta \rho_{-})$ \; (\New{cyclic shift}),
  \item\label{rel-R3}
    $\tau(\rho) = \tau ((\rho_{+}, \rho_{-})) := (\tau \rho_{+}, \tau
    \rho_{-})$ \; (\New{transposition of the first two rows}),
  \item\label{rel-R4}
    $\pr(\rho) = \pr ((\rho_{+}, \rho_{-})) := (\pr \rho_{+}, \pr
    \rho_{-})$ \;(\New{deletion of the first row}),
  \item\label{rel-R5}
    $\rho \times \rho' = (\rho_{+}, \rho_{-}) \times (\rho'_{+},
    \rho'_{-}) := (\rho_{+} \times \rho'_{+}, \rho_{-} \times
    \rho'_{-})$ \; (\New{Cartesian product}),
  \item\label{rel-R6}
    $\rho \cap \rho' = (\rho_{+}, \rho_{-}) \cap (\rho'_{+},
    \rho'_{-}) := (\rho_{+} \cap \rho'_{+}, \rho_{-} \cap
    \rho'_{-})$ \; (\New{intersection}),

  \item\label{rel-R7}
    $\mu_{-}(\rho) = \mu_{-}((\rho_{+}, \rho_{-})) := (\rho_{-}, \rho_{+})$ \; (\New{index translation by
      `$-$\!'}).
\end{enumerate}

  For $Q \subseteq \pmRel(A)$ the \New{$\pm$\hyp{}relational clone
    generated by $Q$}, i.e., the least $\pm$\hyp{}relational clone
  containing $Q$, is denoted by $\pmRg{Q}$ or $\pmRg[A]{Q}$.
\end{definition}

It is possible to construct from the operations of
Definition~\ref{rel-rules} many other operations under which a
$\pm$\hyp{}relational clone is closed (see Remark~\refI{2.14}). In
particular, we will need in the proofs of later theorems the
\New{projection} of an $m$\hyp{}ary relation $\rho$ to rows
$z_1, \dots, z_t \in \{1, \dots, m\}$, which is the $t$\hyp{}ary relation
\[
  \pr_{z_1, \dots, z_i}(\rho) := \{ \, (a_{z_1}, \dots, a_{z_t}) \in
  A^t \mid (a_1, \dots, a_m) \in \rho \, \}.
\]
This naturally extends to $\pm$\hyp{}relations by componentwise
application: for $\rho = (\rho_{+}, \rho_{-}) \in \pmRela[m](A)$,
$\pr_{z_1, \dots, z_t}(\rho) := (\pr_{z_1, \dots, z_t}(\rho_{+}),
\pr_{z_1, \dots, z_t}(\rho_{-}))$.

\begin{definition}\label{def:preserves}
Let $f \in \pmOpa[n](A)$ with $\sgn(f) = (s_1, \dots, s_n)$ and
$\rho = (\rho_{+}, \rho_{-}) \allowbreak \in \pmRela[m](A)$. We say that $f$
\emph{$\pm$\hyp{}preserves} $\rho$, or $f$ is a \emph{$\pm$\hyp{}polymorphism}
of $\rho$, or $\rho$ is \emph{invariant} for $f$, and we write
$f \pmpreserves \rho$, if
\[
  \forall s \in \{ \mathord{+}, \mathord{-} \} : f(\rho_{s_1 s},
  \dots, \rho_{s_n s}) \subseteq \rho_s
\]
(see Definition~\refI{3.4}). In other words, if
$\sgn(f) = (\mathord{+}, \mathord{-}, \dots)$, then
$f \pmpreserves \rho$ if and only if
\[
  f(\rho_{+}, \rho_{-}, \dots ) \subseteq \rho_{+} \qquad \text{and}
  \qquad f(\rho_{-}, \rho_{+}, \dots ) \subseteq \rho_{-}.
\]
\end{definition}

The $\pm$\hyp{}preservation relation induces a Galois connection between
$\pm$\hyp{}operations $\pmOp(A)$ and $\pm$\hyp{}relations $\pmRel(A)$, and the
corresponding operators are $\pmPol Q$ ($\pm$\hyp{}polymorphisms of
$F \subseteq \pmRel(A)$) and $\pmInv F$ (invariant $\pm$\hyp{}relations of
$F \subseteq \pmOp(A)$) (see Definition~\refI{3.6}). It was shown in
Part~I (see Theorems \refI{4.1} and \refI{4.6}) that if $A$ is finite,
then the Galois closures coincide with the $\pm$\hyp{}preclones and the
$\pm$\hyp{}relational clones: for all $F \subseteq \pmOp(A)$ and for all
$Q \subseteq \pmRel(A)$,
\begin{align*}
  \pmSg{F} = \pmPol \pmInv F, \\
  \pmRg{Q} = \pmInv \pmPol Q.
\end{align*}
By the properties of the Galois connection $\pmPol$--$\pmInv$, the
lattice $\pmL_{A}$ of $\pm$\hyp{}preclones is dually isomorphic to the
lattice of $\pm$\hyp{}relational clones. Projections are those
$\pm$\hyp{}operations that $\pm$\hyp{}preserve every $\pm$\hyp{}relation, and
$\pm$\hyp{}diagonals are those $\pm$\hyp{}relations which are invariant for
every $\pm$\hyp{}operation. (See Section~\refI{3}.)

Let $F$ be a $\pm$\hyp{}preclone. The \New{positive part} $F_{+}$
($=F^{[\mathord{+}]}$) of $F$ is the set of all $f \in F$ with
$\sgn(f) = (+, \dots, +)$, the \New{negative part} $F_{-}$
($=F^{[\mathord{-}]}$) of $F$ is the set of all $f \in F$ with
$\sgn(f) = (-, \dots, -)$, and the \New{mixed part} $F_{\pm}$ of $F$
is $F_{\pm} := F \setminus (F_{+} \cup F_{-})$.

The following is a consequence of Propositions~\refI{3.10} and
\refI{4.7}.

\begin{proposition}\label{F1}
  Let $Q \subseteq \pmRel(A)$. Then the positive part $F_{+}$ of the
  $\pm$\hyp{}preclone $F := \pmPol Q$ is a clone given by
  \begin{align*}
    F_{+} = \Pol \{ \sigma \in \Rel(A) \mid
    \exists \sigma' \colon (\sigma,\sigma') \in Q \text{ or } (\sigma',\sigma) \in Q \}.
  \end{align*}
  The clone $\Sg{\mathring{F}}$ generated by the functions
  $f \in \mathring{F}$ \textup{(}ignoring signa\textup{)} is given by
  \begin{align*}
    \Sg{\mathring{F}} = \Pol\{\sigma \in \Rel(A) \mid (\sigma,\sigma) \in \pmRg{Q}\}.
  \end{align*}
\end{proposition}

\begin{lemma}\label{Max-L4}
  If $f \pmpreserves (\sigma,\emptyset)$
  \textup{(}$\sigma \neq \emptyset$\textup{)} then
  $\mathord{+} \in \Sgn(f)$ \textup{(}$f$ cannot have pure negative
  signum\textup{)}.
\end{lemma}

\begin{proof} If
  $\sgn(f) = (\mathord{-}, \dots, \mathord{-}) = (s, \dots, s)$ for
  $s = \mathord{-}$, then $f \pmpreserves \rho = (\sigma,\emptyset)$
  implies
  $\emptyset \neq f(\sigma, \dots, \sigma) = f(\rho_{ss}, \dots,
  \rho_{ss}) \subseteq \rho_{s} = \emptyset$, a contradiction.
\end{proof}

From now on, we will further specialize to the case where
$A = \{0,1\}$. In this special case, $\pm$\hyp{}preclones on $A$ are called
\New{Boolean $\pm$\hyp{}preclones}. The lattice of Boolean $\pm$\hyp{}preclones
is denoted by $\pmL_{2}$.

Since there is only one nontrivial permutation on $A = \{0,1\}$, here
denoted by $\partial \colon 0 \mapsto 1 \mapsto 0$, we speak of
\New{dual} (functions or relations, cf.\ Definition~\refI{5.6})
instead of $\partial$-dual. Recall
$f^{\partial}(x_{1},\dots,x_{n}):=(f(x_{1}^{\partial},\dots,x_{n}^{\partial}))^{\partial}$
with $x^{\partial}:=\partial(x)$. Obviously we have (cf.\ Proposition~\refI{5.10}):

\begin{lemma}\label{lem:partial}
  A $\pm$\hyp{}operation $f$ generates a minimal \textup{(}or maximal, resp.\textup{)}\
  $\pm$\hyp{}preclone if and only if $f^{\partial}$ does.
\end{lemma}


\section{Maximal Boolean $\pm$\hyp{}preclones}\label{sec:Max}

In Proposition~\refI{5.3} conditions for maximal $S$\hyp{}preclones are
given. We specialize these conditions for the Boolean case $A=\{0,1\}$
and $S=\{e,s_{0}\}:=\{+,-\}$ in the following proposition.
\begin{proposition}\label{Max-L1}
  Each maximal $\pm$\hyp{}preclone $F\le\pmOp(A)$ can be characterized as
  $F=\pmPol \rho$ for some $\rho=(\rho_{+},\rho_{-})\in\pmRel(A)$ such
  that we have $(^{*})_{+}$ and $(^{*})_{-}$ defined as follows \textup{(}for
  $s\in\{+,-\}$\textup{)}:
  \begin{itemize}
  \item[$(^{*})_{s}:$] $\rho_{s}$ is a diagonal relation, or \\
    $\Pol\rho_{s}$ is a maximal clone and $\ar(\rho)=\ar(\rho_{s})$ is
    the minimal arity of a nontrivial relation in
    $[\rho_{s}]=\Inv\Pol\rho_{s}$.
  \end{itemize}
  
\end{proposition}

\begin{remark}\label{Max-remark}
  We recall the well\hyp{}known result that there are five maximal Boolean
  clones (consisting of $0$\hyp{}preserving, $1$\hyp{}preserving, monotone,
  self\hyp{}dual and linear operations, respectively) characterizable as
  $\Pol\sigma_{i}$ ($i\in\{0,1,2,3,4\}$) for the following relations:
  \begin{itemize}
  \item $\sigma_{0}:=\{0\}$, $\sigma_{1}:=\{1\}$ (unary),
  \item $\sigma_{2}:=\,\le\,:=\{(x,y)\in A^{2}\mid x\leq y\}$,
    $\sigma_{3}:=\{(0,1),(1,0)\}$ (binary),
  \item $\sigma_{4}:=\{(x,y,z,u)\in A^{4}\mid x+y+z+u=0\}$
    (quaternary).
  \end{itemize}
  These relations can be used for $\rho_{+}$ and $\rho_{-}$ in
  Proposition \ref{Max-L1}, in particular they have minimal arity.
\end{remark}

In Theorem~\ref{Max-theorem} we shall list all $\pm$\hyp{}preclones
from Proposition~\ref{Max-L1} which are really maximal. In preparation we need the
following results.

\begin{lemma}\label{Max-L2}
  $\pmPol(\sigma,\sigma)$ is maximal in $\pmL_{A}$ if \textup{(}and only if\textup{)}
  $\Pol\sigma$ is maximal in $\cL_{A}$.
\end{lemma}

\begin{proof} From Proposition~\refI{5.5} we conclude that the mapping
  $\Phi:\cL_{2}\to\pmL_{2}$ with
  $\Phi(F)=\{f\in\pmOp(A)\mid \mathring f\in F\}$ is an embedding onto
  a filter of $\pmL_{2}$ (here the fact is needed that the monoid
  $S=\{+,-\}$ is a group). Therefore $\Phi(\Pol\sigma)$ is a maximal
  $\pm$\hyp{}preclone whenever $\Pol\sigma$ is maximal. We have
  $\Phi(\Pol\sigma)=\{f\in\pmOp(A)\mid \mathring
  f\preserves\sigma\}=\{f\in\pmOp(A)\mid f\pmpreserves
  (\sigma,\sigma)\}=\pmPol(\sigma,\sigma)$, which completes the
  proof. Note that
  $\mathring f\preserves\sigma\iff f\pmpreserves(\sigma,\sigma)$
  follows from Lemma~\refI{3.5(i)}.
\end{proof}

Recall (Definition \refI{2.2}) that $\neg^{-}\in\pmOpa[1](A)$ is the unary Boolean
function $\neg:0\mapsto 1, 1\mapsto 0$ (negation) with signum
$(-)$. For a relation $\sigma$, we use the notation $\neg(\sigma)$ for
the application of the operation $\neg$ to all tuples in $\sigma$ (and
not the complement of $\sigma$ what in other contexts sometimes is
denoted by $\neg\sigma$). Here $\neg(\sigma)$ coincides with the dual
relation $\sigma^{\partial}$, cf.~Section \ref{sec:prel}, for the
permutation $\partial:A\to A: 0\mapsto 1\mapsto 0$).

\begin{lemma}\label{Max-L3}
  Let $\neg^{-}\pmpreserves\rho=(\sigma,\sigma')$. Then we have:
  \begin{enumerate}[label=\textup{(\alph*)}]
  \item\label{Max-L3a} $\sigma'=\neg(\sigma)$ and
    $\pmRg{\rho}=\{(\alpha,\neg(\alpha))\mid
    \alpha\in\Rg[A]{\{\sigma,\sigma'\}}\}$.
  \item\label{Max-L3b} $\pmPol\rho$ is a maximal $\pm$\hyp{}preclone in
    $\pmL_{2}$ if and only if $\Pol\sigma$ is a maximal clone in
    $\cL_{2}$.
  \end{enumerate}
\end{lemma}
\begin{proof}
  \ref{Max-L3a}: By definition of $\pm$\hyp{}preservation we have
  $\neg(\sigma')\subseteq\sigma$ and $\neg(\sigma)\subseteq\sigma'$,
  therefore
  $\sigma'=\neg\neg(\sigma')\subseteq\neg(\sigma)\subseteq\sigma'$. Consequently
  $\alpha'=\neg(\alpha)$ for every $(\alpha,\alpha')\in\pmRg{\rho}$
  (since $\neg^{-}$ also must $\pm$\hyp{}preserve
  $(\alpha,\alpha')$). Further, $(\alpha,\alpha')\in \pmRg{\rho}$
  implies $\alpha\in\Rg[A]{\sigma,\sigma'}$ (by induction on the
  generating rules from Definition \ref{rel-rules}, note
  $(\alpha',\alpha)\in\pmRg{(\alpha,\alpha')}$ by
  Definition \ref{rel-rules}\ref{rel-R7}). Conversely, if
  $\alpha\in\Rg[A]{\sigma,\sigma'}$ then there must be a pp\hyp{}formula
  $\phi$, such that $\alpha=t_{\phi}(\sigma,\sigma')$ (for notation,
  see Remark~\refI{2.14(C)}). For $\alpha':=t_{\phi}(\sigma',\sigma)$ we get
  $(\alpha,\alpha')=t_{\phi}((\sigma,\sigma'),(\sigma',\sigma))\in\pmRg{(\alpha,\alpha')}$. This
  shows the second part of \ref{Max-L3a}.

  \ref{Max-L3b}: If $\sigma$ is a diagonal $d$ then
  $\rho=(d,\neg(d))=(d,d)$ would be trivial. Therefore $\sigma$ is
  non\hyp{}trivial and the forward implication follows from Proposition \ref{Max-L1}.

  Concerning the backward implication, let $\Pol\sigma$ be a maximal
  clone. Then the clone $\Pol\sigma'$ (with $\sigma'=\neg(\sigma)$) is
  maximal as well. Since $\cL_{2}$ is well known (in particular the
  maximal clones), one can check that $\Pol\sigma$ and $\Pol\sigma'$
  are the only maximal clones containing $\Pol\{\sigma,\sigma'\}$, or
  -- from point of view of relations --, $\Rg{\sigma}$ and
  $\Rg{\sigma'}$ are the only minimal relational clones contained in
  $\Rg{\sigma,\sigma'}$.

  Let $(\alpha,\alpha')\in\pmRg{\rho}$ be nontrivial (we know
  $\alpha'=\neg(\alpha)$ by \ref{Max-L3a}). We show $\rho\in\pmRg{(\alpha,\alpha')}$,
  which gives minimality of $\pmRg{\rho}$ and therefore maximality of
  $\pmPol\rho$. By \ref{Max-L3a} we have
  $\alpha\in\Rg[A]{\sigma,\sigma'}$. Thus at least one of the minimal
  relational clones $\Rg{\sigma}$ or $\Rg{\sigma'}$ must be contained
  in $\Rg{\alpha}$; consequently, there exists $\sigma''$ with
  $(\sigma,\sigma'')\in\pmRg{(\alpha,\alpha')}$ or
  $(\sigma',\sigma'')\in\pmRg{(\alpha,\alpha')}$, respectively (as above, use
  the same pp\hyp{}formula applied to $\alpha'$ as to $\alpha$ in
  order to get $\sigma''$). Since
  $\neg^{-}\pmpreserves (\alpha,\alpha')$, by \ref{Max-L3a} we have
  $\sigma''=\lnot(\sigma)=\sigma'$ or
  $\sigma''=\lnot(\sigma')=\sigma$, respectively. Consequently 
  $\rho\in\pmRg{(\alpha,\alpha')}$ and we are done.
\end{proof}

\begin{theorem}\label{Max-theorem}
  There are nine maximal Boolean $\pm$\hyp{}preclones, as listed below. Each
  such preclone is of the form $F=\pmPol\rho$ for some $\pm$\hyp{}relation
  $\rho=(\rho_{+},\rho_{-})$:

  \begin{enumerate}[label=\textup{(\alph*)}]
  \item\label{G1a} $\pmPol(\sigma,\sigma)$ with
    $\sigma\in\{\sigma_{0},\sigma_{1},\sigma_{2},\sigma_{3},\sigma_{4}\}$,
    cf.\ Remark~\ref{Max-remark}.

    {\rm Note
      $\pmPol(\sigma,\sigma)=\{f^{\lambda}\in \pmOp(A)\mid
      f\in\Pola[n]\sigma, \lambda\in S^{n},n\in \N_{+}\}
      =\{h\in\pmOp(A)\mid \mathring h\in\Pol\sigma\}$.}

  \item\label{G1b} $\pmPol(\leq,\geq)$ \textup{(}all functions where each
    argument with positive signum is order\hyp{}preserving and each
    argument with negative signum is order\hyp{}reversing\textup{)}, cf.\
    Example~\refI{3.7}.

  \item\label{G1c}
    $\pmPol(A,\emptyset)=\pmOp(A)_{+}\cup
    \pmOp(A)_{\pm}=\{f\in\pmOp(A)\mid +\in\Sgn(f)\}$ \textup{(}all functions
    with positive or mixed signum\textup{)}.

  \item\label{G1d} $\pmPol(\nabla_{A},\Delta_{A})$ \textup{(}consists of all
    Boolean $\pm$\hyp{}functions, where each negative argument is fictitious
    \textup{(}including all negative constants\textup{))}, $\nabla_{A}=A^{2}$,
    $\Delta_{A}=\{(x,x)\mid x\in A\}$.
  \item\label{G1e}\label{G1e} $\pmPol(\{0\},\{1\})$.

  \end{enumerate}

\end{theorem}
\begin{proof}

  At first we show that all listed $\pm$\hyp{}preclones are really
  maximal.
 
  \ref{G1a}: By Lemma~\ref{Max-L2} all $\pm$\hyp{}preclones in \ref{G1a}
  are maximal, since $\Pol\sigma_{i}$ ($i\in\{0,\dots,4\}$) are the
  well-known maximal Boolean clones in $\cL_{A}$
  (cf.\ Remark~\ref{Max-remark}).

  \ref{G1b}:
Since $\neg^{-}\pmpreserves(\leq,\geq)$ and
  $\neg(\le)=\; \ge$, it follows from Lemma~\ref{Max-L3}\ref{Max-L3b} that
  $\pmPol(\leq,\geq)$ is maximal. However, maximality of
  $\pmPol(\leq,\geq)$ will be shown also in the
  next Theorem~\ref{thm:Max-L5} in a more general context.

  \ref{G1c}: The application of the ($\pm$\hyp{}relational clone closure)
  rules \ref{rel-R2}--\ref{rel-R7} in Definition~\ref{rel-rules} to $(A,\emptyset)$
  shows that each non\hyp{}diagonal $\pm$\hyp{}relation $(\sigma,\sigma')$ in
  $\pmRg{(A,\emptyset)}$ is of the form $(d,\emptyset)$ or
  $(\emptyset,d)$ for some non\hyp{}empty diagonal $d\in D_{A}$. Obviously
  $(A,\emptyset)\in\pmRg{(d,\emptyset)}$, thus $\pmRg{(A,\emptyset)}$
  is a minimal $\pm$\hyp{}relational clone and therefore
  $\pmPol(A,\emptyset)$ is maximal.

  \ref{G1d}: It is easy to see that $\pmRg{(\nabla_{A},\Delta_{A})}$
  consists of $(\emptyset,\emptyset)$ and all pairs $(d,d')$ of
  non\hyp{}empty diagonals $d,d'\in D_{A}$. As shown below in Case~3 we
  have $(\nabla_{A},\Delta_{A})\in \pmRg{(d,d')}$. Thus
  $\pmRg{(\nabla_{A},\Delta_{A})}$ is minimal and therefore
  $\pmPol(\nabla_{A},\Delta_{A})$ is maximal.

  \ref{G1e}: Since $\neg^{-}\pmpreserves(\{0\},\{1\})$, it follows
  from Lemma \ref{Max-L3}\ref{Max-L3b} that $\pmPol(\{0\},\{1\})$ is
  maximal. 

  Now we show that the list is complete, i.e., there are no more
  maximal clones. Let $F=\pmPol\rho $ be a maximal $\pm$\hyp{}preclone,
  $\rho=(\sigma,\sigma')$. We distinguish several cases according to
  the possible candidates characterized in Proposition~\ref{Max-L1}:
  either the candidate is in the list, i.e., one of the clones
  \ref{G1a}--\ref{G1e}, or we find some larger clone (e.g., by
  indicating a $\pm$\hyp{}relation $\rho'$ derivable from $\rho$, i.e.,
  $\rho'\in\pmRg{\rho}$, and therefore
  $\pmPol\rho\subseteq\pmPol\rho'$).

  \textsl{Case 1:} $\sigma=\sigma'$.

  Because of Lemma~\ref{Max-L2}, $\pmPol(\sigma,\sigma)$ is maximal
  only if $\Pol\sigma$ is a maximal clone. Thus we get the
  $\pm$\hyp{}preclones in \ref{G1a}.

  \textsl{Case 2:} $\sigma\neq\sigma'$, $\sigma'=\emptyset$.

  Then $\sigma$ is not empty and we have $F_{-}=\emptyset$ (no
  negative part) by Lemma~\ref{Max-L4}. Obviously the largest such $F$
  contains \textsl{all} $\pm$\hyp{}functions with positive or mixed
  signa. This is the $\pm$\hyp{}preclone in \ref{G1c}. It is trivial to
  check that each $\pm$\hyp{}function with positive or mixed signum
  $\pm$\hyp{}preserves $(A,\emptyset)$.

  Thus from now on we can assume
  $\emptyset\neq\sigma\neq\sigma'\neq\emptyset$.

  We have to look through the list of all remaining possible
  candidates $\rho=(\sigma,\sigma')$ according to
  Proposition~\ref{Max-L1} (note that
  $\pmPol(\sigma,\sigma')=\pmPol(\sigma',\sigma)$ by
  \ref{rel-rules}\ref{rel-R7}, so only one of the pairs,
  $(\sigma,\sigma')$ or $(\sigma',\sigma)$, has to be considered):

  \textsl{Case 3:} $\rho=(d,d')$ where both $d$ and $d'$ are non\hyp{}empty
  diagonal relations.

  For $m:=\ar(\rho)=1$ we would get the diagonal $\rho=(A,A)$. Thus
  $m\ge 2$. Since $d\neq d'$ there must exist two arguments
  $i,j\in\{1,\dots,m\}$, such that $\pr_{i,j}(d)\ne \pr_{i,j}(d')$,
  therefore
  $\{\pr_{i,j}(d),\pr_{i,j}(d')\}=\{\nabla_{A},\Delta_{A}\}$, i.e.,
  $\pmPol\rho$ is contained in the maximal $\pm$\hyp{}preclone
  $\pmPol(\nabla_{A},\Delta_{A})$ in \ref{G1d}.

  It remains to examine the possible candidates
  $\rho=(\sigma,\sigma')$ where $\sigma$ is nontrivial, and therefore
  $\Pol\sigma$ is a maximal clone by Proposition~\ref{Max-L1}.

  \textsl{Case 4:} $\ar(\rho)=1$ (unary relations)

  There are two unary relations determining a maximal clone,
  $\sigma_{0}$ and $\sigma_{1}$, and one non\hyp{}empty unary trivial
  relation, namely $A$. Thus it remains to consider the following
  pairs:

  \textsl{Case 4a:} $\rho=(\sigma_{0},\sigma_{1})$.

  This gives the maximal $\pm$\hyp{}preclone \ref{G1e}.

  \textsl{Case 4b:} $\rho=(\sigma_{0},A)$ or
  $\rho=(\sigma_{1},A)$. Then $\pmPol\rho$ is contained in the maximal
  $\pm$\hyp{}preclone $(\sigma_{i},\sigma_{i})$ in \ref{G1a}, since
  $\rho \cap \mu_{-}(\rho) = (\sigma_{i},A)\cap
  (A,\sigma_{i})=(\sigma_{i},\sigma_{i})$ ($i\in\{0,1\}$).

  \textsl{Case 5:} $\ar(\rho)=2$ (binary relations)

  There are three binary relations determining a maximal clone,
  $\sigma_{2}$, $\sigma_{2}^{-1}$ and $\sigma_{3}$, and two non\hyp{}empty
  trivial relations, $\Delta_{A}=\{(x,x)\mid x\in A\}$ and
  $\nabla_{A}=A^{2}$. Thus it remains to consider the following pairs:

  \textsl{Case 5a:} $\rho=(\sigma_{2},\sigma_{2}^{-1})$.

  This gives the maximal $\pm$\hyp{}preclone \ref{G1b}.

  \textsl{Case 5b:} $\rho=(\sigma_{2},\Delta_{A})$ (equivalently
  $(\sigma_{2}^{-1},\Delta_{A})$), or $\rho=(\sigma_{3},\Delta_{A})$.

Then $\pmPol\rho$ is contained in the maximal
  $\pm$\hyp{}preclone $\pmPol(\sigma_{i},\sigma_{i})$ in \ref{G1a} since
$\rho\circ\mu_{-}(\rho)=(\sigma_{i},\sigma_{i})$, $i\in\{2,3\}$.

  \textsl{Case 5c:} $\rho=(\sigma_{2},\nabla_{A})$ (equivalently
  $(\sigma_{2}^{-1},\nabla_{A})$) or  $\rho=(\sigma_{3},\nabla_{A})$.

  As in case~4b, $\pmPol\rho$ is contained in the maximal
  $\pm$\hyp{}preclone $\pmPol(\sigma_{i},\sigma_{i})$ since
  $\rho\cap\mu_{-}(\rho)=(\sigma_{i},\nabla_{A})\cap
  (\Delta_{A},\sigma_{i})=(\sigma_{i},\sigma_{i})$, $i\in\{2,3\}$.

\textsl{Case 5d:} $\rho=(\sigma_{2},\sigma_{3})$.

 This case can be reduced to Case~5b since
 $\rho\circ\rho=(\sigma_{2},\Delta_{A})$.

  \textsl{Case 6:} $\ar(\rho)=4$ (quaternary relations).

  There is only one quaternary relation determining a maximal clone,
  namely $\sigma_{4}$, and we have the non\hyp{}empty trivial relations
  $d\in D_{A}^{(4)}$. Thus it remains to consider the following pairs:

  \textsl{Case 6a:} $\rho=(\sigma_{4},A^{4})$.

  As in cases~4b and 5c, $\pmPol\rho$ is contained in the maximal
  $\pm$\hyp{}preclone $\pmPol(\sigma_{4},\sigma_{4})$.

  \textsl{Case 6b:} $\rho=(\sigma_{4},d)$ where $d$ is a quaternary
  diagonal such that there exist two arguments $i,j\in\{1,2,3,4\}$
  such that $\pr_{i,j}(d)=\Delta_{A}$.

  Then $\pr_{i,j}(\rho)=(\nabla_{A},\Delta_{A})$ (since
  $\pr_{i,j}(\sigma_{4})=\nabla_{A}$), i.e., $\pmPol\rho$ is contained
  in the maximal $\pm$\hyp{}preclone $\pmPol(\nabla_{A},\Delta_{A})$ in
  \ref{G1d}. Consequently $\pmPol\rho$ is not maximal.
\end{proof}

As already mentioned, part \ref{G1b} of Theorem~\ref{Max-theorem}
is a special case of the following theorem where we shall see that our
``motivating example'', Example \refI{2.5}, is maximal not only for Boolean
$\pm$\hyp{}preclones but also in general for an arbitrary finite base set $A$.

\begin{theorem}\label{thm:Max-L5}
    Let $(A,\leq)$ be a finite poset with 0 and 1 and $S=\{+,-\}$. Then
  $\pmPol(\leq,\geq)$ is a maximal $\pm$\hyp{}preclone in $\pmL_{A}$.
\end{theorem}

\begin{proof} This proof is a generalization of the proof of the
  corresponding result for usual clones, see
  \cite[Satz~4.3.7]{PoeK79}.

Let $\rho=(\rho_{+},\rho_{-}):=(\leq,\geq)$, $F:=\pmPol\rho$ and let
$\sigma=(\sigma_{+},\sigma_{-})\in\pmInva[m] F\setminus\pmD_{A}$ be a
non\hyp{}trivial (\refI{2.9}) $m$\hyp{}ary $\pm$\hyp{}relation. We are going to show
$\rho\in\pmRg[A]{\sigma}$, which will prove that $\pmRg[A]{\rho}$ is a
minimal $\pm$\hyp{}relational clone; therefore $\pmPol\rho$ is maximal and we are
done.

Without loss of generality, we may assume that
$\sigma$ has no equal rows (otherwise one could delete one of the
rows and proceed with the resulting $\pm$\hyp{}relation), i.e., we assume 
$\pr_{i,j}(\sigma)\not\subseteq
(\Delta_{A},\Delta_{A})$ for all $i,j\in\{1,\dots,m\}$.

At first we need some notation for $\sigma$ considered as a matrix. 
Let $\sigma_{+} = \allowbreak \{r_{1}, \allowbreak \dots, \allowbreak r_{n}\}$,
$\sigma_{-}=\{r'_{1},\dots,r'_{n'}\}$ (i.e., we consider the colums
$r_{i}=(r_{i}(1),\dots,r_{i}(m))\in A^{m}$ and $r'_{i}\in A^{m}$ of $\rho$) and let 
$(\bz_{i},\bz'_{i}):=((r_{1}(i),\dots,r_{n}(i)),(r'_{1}(i),\dots,r'_{n'}(i)))$ denote the
rows of $\sigma$, $i\in\{1,\dots,m\}$. Moreover, writing 
$\bz_{i}\leq \bz_{j}$ is understood componentwise, i.e., $r_{k}(i)\leq r_{k}(j)$ for
each $k\in\{1,\dots,n\}$, analogously
$\bz'_{i}\geq\bz'_{j}$ means $r'_{k}(i)\geq r'_{k}(j)$ for each
$k\in\{1,\dots,n'\}$.

We distinguish two cases.

\textsl{Case 1}: {\it There exist $i,j\in\{1,\dots,m\}$, $i\neq j$,
  such that $\bz_{i}\leq\bz_{j}$ and $\bz'_{i}\geq \bz'_{j}$
\textup{(}equivalently, $\pr_{i,j}(\sigma)\subseteq\rho$\textup{)}.}

We are going to show $\sigma':=\pr_{i,j}(\sigma)=\rho$, thus
$\rho\in\pmRg[A]\sigma$ what was to be proved.

In fact, since by the above assumption $\sigma'$ is non\hyp{}trivial, there
exists $(a,b)\in\sigma'_{+}\cup\sigma'_{-}$ such that $a\neq b$. Assume at first that
$(a,b)\in\sigma'_{+}$, thus $a<b$ (since $\bz_{i}\leq\bz_{j}$). Then
for arbitrary $(c,d)\in A^{2}$ with $c\leq d$, the $\pm$\hyp{}operation 
\begin{align*}
  f(x):=\begin{cases}
    c &\text{if }x\leq a,\\
    d &\text{otherwise}
  \end{cases}
\end{align*}
with $\sgn(f)=(+)$
belongs to $F$ (it is easy to check that it $\pm$\hyp{}preserves $\rho$)
and therefore $(c,d)=(f(a),f(b))\in\sigma'_{+}$.
Likewise, the  $\pm$\hyp{}operation 
\begin{align*}
  g(x):=\begin{cases}
    c &\text{if }x\geq b,\\
    d &\text{otherwise}
  \end{cases}
\end{align*}
with $\sgn(g)=(-)$ also belongs to $F$ and we get
$(d,c)=(g(a),g(b))\in\sigma'_{-}$. Since $c\leq d$ was chosen arbitrarily, we
have
$\rho=(\leq,\geq)\subseteq(\sigma'_{+},\sigma'_{-})=\sigma'$. Thus
$\rho=\sigma'$ since $\sigma'\subseteq\rho$ by assumption of Case 1,
and we are done. 

In the remaining case that $(a,b)\in\sigma'_{-}$ (i.e., $a>b$) we can
apply the above proof to $\mu_{-}(\sigma)$ and
$\sigma'':=\pr_{j,i}(\mu_{-}(\sigma))$ yielding
$\rho=\sigma''$, consequently we also have $\rho=\sigma''\in\pmRg[A]\sigma$.

\textsl{Case 2}: {\it All rows are ``incomparable'' in the sense
  that $\bz_{i}\leq \bz_{j}$ and $\bz'_{i}\geq \bz'_{j}$ implies $i=j$
  \textup{(}for $i,j\in\{1,\dots,m\}$\textup{)}.}

We shall show that $\sigma={(A^{m},A^{m})}$,
i.e., $\sigma$ is trivial, which contradicts the assumption
about $\sigma$. Thus this case cannot appear.

Let $(b_{1},\dots,b_{m})\in A^{m}$. Consider the $(n+n')$\hyp{}ary
$\pm$\hyp{}operation $f$ with signum $\sgn(f)=(+,\ldots,+,-,\ldots,-)$
(with $n$ signs `$+$' and $n'$ signs `$-$') defined by
\begin{align*}
  f(\bx,\bx'):=
  \begin{cases}
    b_{i}&\text{if }(\bx,\bx')=(\bz_{i},\bz'_{i}) 
      \text{ for } i\in\{1,\dots,m\},\\
    0 &\text{if }\exists\,i\in\{1,\dots,m\}:
    \bx\leq\bz_{i},\bx'\geq\bz'_{i}\text{ and }(\bx,\bx')\neq(\bz_{i},\bz'_{i}),\\ 
    1&\text{otherwise.}
  \end{cases}
\end{align*}
Here $(\bx,\bx')=((x_{1},\dots,x_{n}),(x'_{1},\dots,x'_{n'}))\in
A^{n}\times A^{n'}$. 
Then $f$ is well defined and it $\pm$\hyp{}preserves
$\rho$ (decreasing (with respect to $\leq$) the value of some $x_{k}$
or increasing the value of some $x'_{k}$ cannot increase the
value of $f$) --- for both properties the above ``incomparability'' of the rows
is needed and indispensable.
 Observe that
$(b_{1},\dots,b_{m})=(f(\bz_{1},\bz'_{1}),\dots,f(\bz_{m},\bz'_{m}))$
equals to $f(r_{1},\dots,r_{n},r'_{1},\dots,r'_{n'})$ ($f$ applied
componentwise), but the latter must belong to $\sigma_{+}$ since
$f\pmpreserves\sigma$ and
$r_{1},\dots,r_{n}\in\sigma_{+}$ and
$r'_{1},\dots,r'_{n'}\in\sigma_{-}$.

Analogously one can prove that $(b_{1},\dots,b_{m})\in\sigma_{-}$
using the following $\pm$\hyp{}operation $g$ with signum $(-,\ldots,-,+,\ldots,+)$
(with $n$ signs `$-$' and $n'$ signs `$+$')

\begin{align*}
  g(\by,\by'):=
  \begin{cases}
    b_{i}&\text{if }(\by,\by')=(\bz_{i},\bz'_{i}) 
      \text{ for } i\in\{1,\dots,m\},\\
    0 &\text{if
    }\exists\,i\in\{1,\dots,m\}:\by\geq\bz_{i},\by'\leq\bz'_{i}\text{
      and }(\by,\by')\neq(\bz_{i},\bz'_{i}),\\ 
    1&\text{otherwise.}
  \end{cases}
\end{align*}
Then
$(b_{1},\dots,b_{m})=(g(\bz_{1},\bz'_{1}),\dots,g(\bz_{m},\bz'_{m}))=g(r_{1},\dots,r_{n},r'_{1},\dots,r'_{n'})\in\sigma_{-}$
because $g\pmpreserves\rho$. Since $(b_{1},\dots,b_{m})$ was chosen
arbitrarily we have $\sigma=(A^{m},A^{m})$.
\end{proof}

The full description of maximal $\pm$\hyp{}preclones in
Theorem~\ref{Max-theorem} makes the completeness 
problem tractable. A set $F\subseteq \pmOp(A)$ of $\pm$\hyp{}operations is
called \New{complete} if it generates all $\pm$\hyp{}operations: $\pmSg{F}=\pmOp(A)$. 

\begin{theorem}[Completeness theorem for $\pm$\hyp{}preclones]
\label{thm:complete1} 
A set $F\subseteq\pmOp(A)$ of $\pm$\hyp{}operations
     is complete 
    if and only for each of the following nine $\pm$\hyp{}relations there
    exists a $\pm$\hyp{}operation in $F$ which does not preserve this
    relation: 
    $(\sigma_{i},\sigma_{i})$ for $i\in\{0,1,2,3,4\}$,
    $(\leq,\geq)$, $(A,\emptyset)$,
    $(\nabla_{A},\Delta_{A})$, $(\{0\},\{1\})$.
\end{theorem}

\begin{proof} The lattice $\pmL_{2}$ is coatomic by
  Proposition~\refI{5.1(3)}. Therefore $F$ is complete if and only if $F$
  (and therefore also $\pmSg{F}$) is
  not contained in a maximal $\pm$\hyp{}preclone. Since, by
  Theorem~\ref{Max-theorem}, each maximal $\pm$\hyp{}preclone is of the
  form $M=\pmPol\rho$ 
  for $\rho$ being one of the indicated nine $\pm$\hyp{}relations, it
  enough to find a $\pm$\hyp{}operation $f\in F$ which does not preserve $\rho$.
\end{proof}


\section{$S$\hyp{}algebras and $S$\hyp{}terms}\label{sec:Sterms}

For the proofs of some results in the subsequent sections, we need the
notions of $S$\hyp{}algebras and terms and term operations in languages of
$S$\hyp{}algebras. These are rather straightforward adaptations of the
classical notions of terms and term operations in usual algebras. Some
care is, however, required, because the expected recursive definitions
do not work well for arbitrary monoids $S$.

\begin{definition}\label{def:S-term}

  Let $S$ be a fixed finite monoid. An \New{$S$\hyp{}algebra} is a
  structure $\mathbf{A} = (A, (f_i)_{i \in I})$, where $A$ is a set
  and $(f_i)_{i \in I}$ is a family of $S$\hyp{}operations on $A$, each
  $f_i$ of signum $\lambda_i = (s^i_1, \dots, s^i_{n_i})$. The family
  $(\lambda_i)_{i \in I}$ is the \New{type} of the $S$\hyp{}algebra
  $\mathbf{A}$.

  For a type $\tau = (\lambda_i)_{i \in I}$ and a set $X$ of
  variables, we let $(f_i)_{i \in I}$ be a set of operation symbols,
  and we define \New{terms of type $\tau$ over $X$} in the usual
  recursive way. Each variable $x_i \in X$ is a term, and if $f_i$ is
  an operation symbol with signum
  $\sgn(f_i) = \lambda_i = (s^i_1, \dots, s^i_{n_i})$ and
  $t_1, \dots, t_{n_i}$ are terms, then $f_i(t_1, \dots, t_{n_i})$ is
  a term. (In other words, these are just terms of type
  $(n_i)_{i \in I}$ in the usual universal-algebraic sense.) Here,
  $f_i$ is called the \New{root symbol} and $t_1, \dots, t_{n_i}$ are
  called the \New{immediate subterms} of the term
  $f_i(t_1, \dots, t_{n_i})$. The root symbol of a term $x_i \in X$ is
  $x_i$, and this term has no immediate subterms.
\end{definition}

\begin{definition}[$S$\hyp{}terms]
  For any term $t$ of type $\tau$, we assign a sign to each occurrence
  of a variable in $t$ by the following recursion. If $t = x_i$ is a
  variable, then the single variable $x_i$ gets the sign $e$ (the
  identity element of $S$). If $t = f_i(t_1, \dots, t_n)$, where $f_i$
  is an $S$\hyp{}operation symbol with $\sgn(f_i) = (s^i_1, \dots, s^i_n)$,
  then an occurrence (in $t$) of a variable $x$ that lies within the
  immediate subterm $t_j$ gets the sign $ss^i_j$ (product in $S$),
  where $s$ is the sign of $x$ in the term $t_j$. A term $t$ is called
  an \New{$S$\hyp{}term} if for each variable $x_i \in X$, every occurrence of $x_i$ in $t$ has the
  same sign.
\end{definition}

\begin{definition}[recursive $S$\hyp{}terms and $S$\hyp{}term operations]
  Let $X_{n} = \{x_{1}, \allowbreak \dots, \allowbreak x_{n}\}$ (set of variable symbols). A
  \New{recursive $S$\hyp{}term $t$ over $X_{n}$} with signum
  $(s_1, \dots, s_n)$ is defined recursively as follows:
  \begin{enumerate}
  \item Each variable $x_{i}$ ($i \in \{1, \dots, n\}$) is a recursive
    $S$\hyp{}term over $X_{n}$ with signum
    $(s_1, \dots, s_{i-1}, e, s_{i+1}, \dots, s_n)$, for all
    $s_1, \dots, s_{i-1}, s_{i+1}, \dots, s_n \in S$.
  \item Let $f$ be an $m$\hyp{}ary operation symbol with
    $\sgn(f) = (s_1, \dots, s_m)$ and, for $j \in \{1, \dots, m\}$,
    let $t_{j}$ be a recursive $S$\hyp{}term over $X_{n}$ with signum
    $(s_{j,1},\dots,s_{j,n})$ such that
    $s_{j,i} s_{j} = s_{j',i} s_{j'} =: s'_i$ for all
    $i \in \{1, \dots, n\}$ and $j, j' \in \{1, \dots, m\}$. Then
    $f(t_{1},\dots,t_{m})$ is a recursive $S$\hyp{}term over $X_{n}$ with
    signum $(s'_1, \dots, s'_n)$.
  \end{enumerate}

  A recursive $S$\hyp{}term $t$ over $X_n$ with signum $(s_1, \dots, s_n)$
  induces an \New{$S$\hyp{}term operation} $t^\mathbf{A}$ with signum
  $(s_1, \dots, s_n)$ on an $S$\hyp{}algebra $\mathbf{A}$ via the usual
  recursion.
  \begin{enumerate}
  \item If $t = x_i$ ($i \in \{1, \dots, n\}$), then
    $t^\mathbf{A} = \pr^{(m)}_i$.
  \item If $t = f(t_1, \dots, t_m)$, where the signum of $f$ and those
    of the subterms $t_1, \allowbreak \dots, \allowbreak  t_m$ satisfy the conditions of the
    above definition,
    then
    $t^\mathbf{A} = f^\mathbf{A}(t_1^\mathbf{A}, \dots,
    t_m^\mathbf{A})$.
  \end{enumerate}
\end{definition}

\begin{proposition}\label{prop:S-term}
  \leavevmode
  \begin{enumerate}[label={\upshape(\roman*)}]
  \item\label{prop:S-term:1} Every recursive $S$\hyp{}term of type $\tau$
    over $X_n$ is an $S$\hyp{}term. Moreover, if $t$ has signum
    $(s_1, \dots, s_n)$, then, for $i \in \{1, \dots, n\}$, every
    occurrence of $x_i$ has sign $s_i$ in $t$.
  \item\label{prop:S-term:sub} If $S$ is a right cancellative monoid,
    then every subterm of an $S$\hyp{}term is an $S$\hyp{}term.
  \item\label{prop:S-term:2} If $S$ is a right cancellative monoid,
    then every $S$\hyp{}term of type $\tau$ over $X_n$ is a recursive
    $S$\hyp{}term.
  \end{enumerate}
\end{proposition}

\begin{proof}
  \ref{prop:S-term:1} We prove by induction on the complexity of terms
  that if $t$ is a recursive $S$\hyp{}term of signum $(s_1, \dots, s_n)$,
  then, for $i \in \{1, \dots, n\}$, every occurrence of $x_i$ has
  sign $s_i$ in $t$; the claim follows from this. If $t = x_i$ is a
  recursive $S$\hyp{}term of signum
  $(s_1, \dots, s_{i-1}, e, s_{i+1}, \dots, s_n)$, then clearly $x_i$
  has sign $e$ in $t$. Let now $t = f(t_1, \dots, t_m)$ with
  $\sgn(f) = (s_1, \dots, s_m)$ and, for $j \in \{1, \dots, m\}$,
  $t_j$ is a recursive $S$\hyp{}term with signum
  $(s_{j,1}, \dots, s_{j,n})$ such that $s_{j,i} s_j = s'_i$ for all
  $i \in \{1, \dots, n\}$. By the inductive hypothesis, the sign of
  every occurrence of $x_i$ in $t_j$ is $s_{j,i}$; each occurrence of
  $x_i$ within $t_j$ has sign $s_{j,i} s_j = s'_i$ in $t$ -- this
  holds for all $j \in \{1, \dots, m\}$.

  \ref{prop:S-term:sub} The claim will follow if we show that the
  immediate subterms of an $S$\hyp{}term are $S$\hyp{}terms. Let
  $t = f ( t_1, \dots, t_n )$ be an $S$\hyp{}term, where $f$ is an
  operation symbol with signum $(s_1, \dots, s_n)$. Assume that for
  some $j \in \{1, \dots, n\}$, the immediate subterm $t_j$ contains
  two occurrences of a variable $x$ with signs $s$ and $s'$ (in
  $t_j$). These occurrences of $x$ have signs $s s_j$ and $s' s_j$ in
  $t$. Because $t$ is an $S$\hyp{}term, we have $s s_j = s' s_j$. By the
  right cancellation property we get $s = s'$. Consequently, each
  immediate subterm $t_j$ of $t$ is an $S$\hyp{}term.

  \ref{prop:S-term:2} We prove by induction on the complexity of terms
  that if $t$ is an $S$\hyp{}term and $(s_1, \dots, s_n) \in S^n$ is such
  that, for $i \in \{1, \dots, n\}$, every occurrence of $x_i$ has
  sign $s_i$ in $t$, then $t$ is a recursive $S$\hyp{}term of signum
  $(s_1, \dots, s_n)$; note that such an induction is possible by part
  \ref{prop:S-term:sub}. The basis of the induction is clear:
  $t = x_i$ is an $S$\hyp{}term, where the unique occurrence of $x_i$ has
  sign $e$, and it is a recursive $S$\hyp{}term of every signum
  $(s_1, \dots, s_{i-1}, e, s_{i+1}, \dots, s_n)$. Assume now that
  $t = f(t_1, \dots, t_m)$, where $\sgn(f) = (s_1, \dots, s_m)$ and
  that, for $i \in \{1, \dots, n\}$, every occurrence of $x_i$ has
  sign $s'_i$ in $t$. By part \ref{prop:S-term:sub}, for
  $j \in \{1, \dots, m\}$, $t_j$ is an $S$\hyp{}term, so, for
  $i \in \{1, \dots, n\}$, every occurrence of $x_i$ within $t_j$ has
  sign $s_{j,i}$ in $t_j$ for some $s_{j,i} \in S$ such that
  $s_{j,i} s_j = s'_i$. By the inductive hypothesis, $t_j$ is a
  recursive $S$\hyp{}term of signum $(s_{j,1}, \dots,
  s_{j,n})$. Consequently, $t = f(t_1, \dots, t_m)$ is a recursive
  $S$\hyp{}term of signum $(s'_1, \dots, s'_n)$.
\end{proof}

\begin{remark}\label{rem:Slcm}
  If the monoid $S$ fails to be right cancellative, then subterms of
  an $S$\hyp{}term may fail to be $S$\hyp{}terms. For example, let
  $S = (\{0,1\}, \land)$ be the two\hyp{}element semilattice with neutral
  element $1$, and consider the type of $S$\hyp{}algebras comprising one
  operation symbol $f$ with signum $(0,1)$. Then $f(f(x_1,x_1),x_2)$
  is an $S$\hyp{}term, as both occurrences of $x_1$ have sign $0$ and the
  occurrence of $x_2$ has sign $1$, but its subterm $f(x_1,x_1)$ fails
  to be an $S$\hyp{}term.

  Even in the non\hyp{}right\hyp{}cancellative case, it makes sense to consider
  $S$\hyp{}term operations induced by $S$\hyp{}terms. However, we cannot write a
  recursive definition purely in terms of $S$\hyp{}operations (because
  subterms of $S$\hyp{}terms may fail to be $S$\hyp{}terms). We can nevertheless
  work out the recursion with the underlying operations and assign the
  signum in the end. Namely, if $t$ is an $S$\hyp{}term over $X_n$ and
  $(s_1, \dots, s_n)$ satisfies the condition that, for
  $i \in \{1, \dots, n\}$, every occurrence of $x_i$ has sign $s_i$ in
  $t$, then we can say that $t$ has signum $(s_1, \dots, s_n)$ and
  that the $S$\hyp{}term operation $t^\mathbf{A}$ has signum
  $(s_1, \dots, s_n)$.
\end{remark}

We will be mainly concerned with the case when $S$ is a right
cancellative monoid. By Proposition~\ref{prop:S-term}, all $S$\hyp{}terms
over $X_n$ are recursive in this case.

We leave it as an exercise to the reader to prove that for an
$S$\hyp{}algebra $\mathbf{A} = (A, (f_i^{\mathbf{A}})_{i \in I})$, the set
of all $S$\hyp{}operations induced on $\mathbf{A}$ by $S$\hyp{}terms is an
$S$\hyp{}preclone. Moreover, every $S$\hyp{}preclone arises in this way from
some $S$\hyp{}algebra.

\begin{definition}\label{def:Snearprojection}
  A nontrivial function $f\in\SOp(A)$ is called an
  \New{$S$\hyp{}near\hyp{}projection} if identification of arguments with the
  same sign (cf.\ Definition~\refI{2.3(5)} and Remark~\refI{2.4}) always gives a trivial
  function, i.e., a projection in $\SJ_{A}$, more precisely, if for
  all $p, q \in \{1, \dots, n\}$ such that $p \neq q$ and the $p$-th
  and the $q$-th arguments have the same sign, there exists an
  $i \in \{1, \dots, n\}$ such that the $i$-th argument has sign $e$
  (the identity element of $S$) and $f(x_1, \dots, x_n) = x_i$
  whenever $x_p = x_q$. A nontrivial function $f \in \SOpa[n]_{A}$ is
  called an \New{$S$\hyp{}semiprojection} if there exists
  $i \in \{1,\dots,n\}$ such that the $i$-th argument has sign $e$ and
  $f(x_{1},\dots,x_{n})=x_{i}$ whenever there exist two equal arguments with
  the same sign.
\end{definition}

\begin{lemma}\label{lem:Snearprojection1}
  Let $S$ be a right cancellative monoid with identity element
  $e$. Let $f$ be an $n$\hyp{}ary $S$\hyp{}near\hyp{}projection. Then $\SSg{f}$
  contains no non\hyp{}trivial $S$\hyp{}operations with arity less than $n$.
\end{lemma}

\begin{proof}
  Let $f$ be an $n$\hyp{}ary $S$\hyp{}near\hyp{}projection with signum
  $\lambda = (s_1, \dots, s_n)$. Let $\mathbf{A} = (A, f)$ be an
  $S$\hyp{}algebra of type $\tau = (\lambda)$. The statement follows if we
  prove that for every positive integer $m$ less than $n$, every
  $S$\hyp{}term $t$ of type $\tau$ over $X_m$ induces on $\mathbf{A}$ a
  trivial $S$\hyp{}operation, that is, $t^{\mathbf{A}} = \pr^{(m)}_i$ for
  some $i \in \{1, \dots, m\}$ such that the sign of the essential
  $i$-th argument is $e$. We prove this by induction on the complexity
  of terms.

  If $t = x_i \in X_m$, then $t^{\mathbf{A}} = \pr^{(m)}_i$ and the
  essential $i$-th argument has sign $e$ because the unique variable
  of the term $x_i$ has sign $e$ by definition.

  Consider now the case when $t = f ( t_1, \dots, t_n )$. For
  $j = 1, \dots, n$, the immediate subterm $t_j$ is an $S$\hyp{}term by
  Lemma~\ref{prop:S-term}\ref{prop:S-term:sub}. By the inductive
  hypothesis, for $j = 1, \dots, n$, we have that
  $t_j^{\mathbf{A}} = \pr^{(m)}_{i_j} = x_{i_j}^{\mathbf{A}}$ for some
  $i_j \in \{1, \dots, m\}$ and the variable $x_{i_j}$ has sign $e$ in
  $t_j$. Thus, the occurrences of $x_{i_j}$ in $t$ within the
  immediate subterm $t_j$ have sign $e s_j = s_j$. In fact, every
  occurrence of $x_{i_j}$ in $t$ has sign $s_j$ because $t$ is an
  $S$\hyp{}term. Consequently, the term $f(x_{i_1}, \dots, x_{i_n})$ is an
  $S$\hyp{}term, and we have
  $t^{\mathbf{A}} = (f(x_{i_1}, \dots,
  x_{i_n}))^{\mathbf{A}}$. Because $m < n$, the variables
  $x_{i_1}, \dots, x_{i_n}$ are not all distinct. Since $f$ is an
  $S$\hyp{}near\hyp{}projection, it follows that
  $(f(x_{i_1}, \dots, x_{i_n}))^{\mathbf{A}} = \pr^{(m)}_{i_\ell} =
  x_{i_\ell}^\mathbf{A}$ for some $\ell \in \{1, \dots, n\}$ with
  $s_\ell = e$, so the sign of the essential $i_\ell$-th argument is
  $e$. We conclude that $t^{\mathbf{A}}$ is a trivial $S$\hyp{}operation.
\end{proof}

\begin{remark}
  Lemma~\ref{lem:Snearprojection1} may fail for monoids $S$ that are
  not right cancellative. For a counterexample, consider the
  non\hyp{}right\hyp{}cancellative monoid $(S, \mathord{\circ})$, where
  $S = \{e, a, b\}$ and $e \circ x = x \circ e = x$ for all $x \in S$,
  and $x \circ y = a$ for all $x, y \in \{a, b\}$, and the
  $S$\hyp{}operation $f$ of signum $(a,b)$ on $\{0,1\}$, given by the rule
  $f(x,y) = 0$ for all $x, y \in \{0,1\}$. The $S$\hyp{}operation $f$ is an
  $S$\hyp{}near\hyp{}projection; the condition of Definition
  \ref{def:Snearprojection} is vacuously satisfied. However, $\SSg{f}$
  contains the $S$\hyp{}operation induced by the term $t = f(x,f(x,x))$ on
  $\mathbf{A} = (\{0,1\}, f)$, and this term operation is the unary
  constant 0 operation of signum $(a)$, a non\hyp{}trivial $S$\hyp{}operation
  with arity less than that of $f$.
\end{remark}


\section{Minimal Boolean $\pm$\hyp{}preclones}\label{sec:Min}

 In this section we shall describe all minimal $\pm$\hyp{}preclones, i.e.,
 the atoms of the lattice $\pmL_{A}$. Minimality of a $\pm$\hyp{}preclone
 $F$ is characterized by the property that each non\hyp{}trivial
 $\pm$\hyp{}operation $f \in F$ already generates $F$: $\pmSg{f} = F$. A first
 observation is:

\begin{lemma}\label{lem:min1}
  Each minimal $\pm$\hyp{}preclone $F$ has trivial positive part \textup{(}only
  projections\textup{)} and empty negative part: $F_{+}\subseteq \pmJ_{A}$,
  $F_{-}=\emptyset$.
\end{lemma}

\begin{proof} If $F=\pmSg{f}$ for some nontrivial function $f$ with
  positive signum $\sgn(f)=(+,\dots,+)$ or negative signum
  $\sgn(f)=(-,\dots,-)$, then $\pmSg{f}$ contains a nontrivial proper
  $\pm$\hyp{}preclone, namely $\pmSg{\nabla^{-}f}\subset\pmSg{f}$ or
  $\pmSg{\nabla^{+}f}\subset\pmSg{f}$, respectively (for the notation, see
  Section~\ref{sec:prel} or Definition~\refI{2.3(4)}).
\end{proof}

\begin{definition}[minimal signum]\label{def:nearprojection}

  Recall that $\ar(f)$ denotes the arity of $f$. Let $\essar(f)$ denote the
  number of essential arguments of $f$, called the \New{essential arity} of $f$.

  The \New{minimal signum} of a minimal $\pm$\hyp{}preclone $F$ is the
  signum of a $\pm$\hyp{}operation $f\in F$ with minimal essential arity,
  and with the greatest number of entries equal to $+$. Then we say that $f$ \New{has minimal
  signum} and assume without loss of generality that the entries of the signum are ordered
  such that every $+$ precedes every $-$ (otherwise permute the arguments accordingly).
\end{definition}

\begin{remark}\label{rem:minsignum}
  Each $\pm$\hyp{}operation $f \in F$ with minimal signum must be a
  $\pm$\hyp{}near\hyp{}projection (see Definition \ref{def:Snearprojection},
  specialized to the case $S = \pm$), for otherwise we could further
  identify arguments of equal sign, getting a nontrivial
  $\pm$\hyp{}operation of smaller arity.
\end{remark}

\begin{lemma}\label{lem:nearprojection} Each minimal $\pm$\hyp{}preclone
  $F$ contains a $\pm$\hyp{}near\hyp{}projection $f$ with $2 \leq \ar(f) \leq 4$
  and with one of the following minimal signa:
  $(\mathord{+},\mathord{+},\mathord{+},\mathord{-})$,
  $(\mathord{+},\mathord{+},\mathord{-},\mathord{-})$,
  $(\mathord{+},\mathord{+},\mathord{-})$,
  $(\mathord{+},\mathord{-},\mathord{-})$,
  $(\mathord{+},\mathord{-})$. Moreover, in the case when the signum
  is $(\mathord{+},\mathord{+},\mathord{+},\mathord{-})$, the argument
  with sign $-$ is fictitious.
\end{lemma}

\begin{proof}
Assume that $f \in F$ has arity $n$ and minimal signum $\sgn(f)=(s_{1},\dots,s_{n})$.
By Remark~\ref{rem:minsignum}, $f$ must be a $\pm$\hyp{}near\hyp{}projection.
Note that, by
  Lemma~\ref{lem:min1}, there must be at least one $+$ and at least
  one $-$ in the signum, so $n \geq 2$. Furthermore, note that if
  there is a fictitious argument in $f$ with sign $s$, then it must be
  the only argument with this sign (otherwise we could identify the
  fictitious argument with another argument with the same sign and
  would get a nontrivial function, in contradiction to the minimality
  of the signum).

  \textit{Claim.} If $\sgn(f)$ has three equal signs, then $n = 4$,
  $\sgn(f) = (\mathord{+}, \mathord{+}, \mathord{+}, \mathord{-})$,
  and the argument with sign $-$ is fictitious.

  \textit{Proof of the Claim.} Assume that $\sgn(f)$ has three equal
  signs, say $s_p = s_q = s_r$ for some distinct
  $p, q, r \in \{1, \dots, n\}$. Because $f$ is a
  $\pm$\hyp{}near\hyp{}projection, there are indices $i$, $j$, and $k$ such that
  $s_i = s_j = s_k = \mathord{+}$ and
  \begin{equation}
    f(x_1, \dots, x_n) =
    \begin{cases}
      x_i & \text{whenever $x_p = x_q$,} \\
      x_j & \text{whenever $x_p = x_r$,} \\
      x_k & \text{whenever $x_q = x_r$.}
    \end{cases}
    \label{eq:f-pqr}
  \end{equation}
  Because $\lvert A \rvert = 2$, there are always two equal arguments
  among $x_p$, $x_q$, and $x_r$. Therefore, we see from
  \eqref{eq:f-pqr} that $f$ may depend only on the arguments indexed
  by $p, q, r, i, j, k$; all other arguments are fictitious.

  Consider first the case that one of $i$, $j$, $k$ is not in
  $\{p, q, r\}$ (note that this is necessarily the case when
  $s_p = s_q = s_r = \mathord{-}$), then by considering tuples
  $(x_1, \dots, x_n)$ with $x_p = x_q = x_r$, we can conclude that
  $i = j = k$. Consequently, $f$ is the $i$-th projection, which
  contradicts the nontriviality of $f$.

  Consider now the case that $i, j, k \in \{p, q, r\}$. Then
  $s_p = s_q = s_r = \mathord{+}$ and the only essential arguments of
  $f$ are the ones indexed by $p, q, r$. Because $\sgn(f)$ contains
  both signs $+$ and $-$ and because a fictitious argument must be the
  only argument of its sign, it follows that $n = 4$ and the
  remaining, fictitious argument is of sign $-$, that is,
  $\sgn(f) = (\mathord{+}, \mathord{+}, \mathord{+},
  \mathord{-})$. This completes the proof of the claim. $\diamond$

  Now, if $n \geq 5$, then there must be three equal signs in
  $\sgn(f)$ by the strong form of the pigeonhole principle, but then
  $n = 4$ by the above Claim, a contradiction. Therefore, $n \leq
  4$. Because the signum
  $(\mathord{+}, \mathord{-}, \mathord{-}, \mathord{-})$ is not
  minimal by the Claim, there remain only the following signa as
  candidates for a minimal signum:
  $(\mathord{+}, \mathord{+}, \mathord{+}, \mathord{-})$,
  $(\mathord{+}, \mathord{+}, \mathord{-}, \mathord{-})$,
  $(\mathord{+}, \mathord{+}, \mathord{-})$,
  $(\mathord{+}, \mathord{-}, \mathord{-})$,
  $(\mathord{+}, \mathord{-})$. Moreover, in the case when the signum
  is $(\mathord{+}, \mathord{+}, \mathord{+}, \mathord{-})$, the
  argument with sign $-$ is fictitious by the Claim.
\end{proof}

\begin{lemma}\label{min-4aryA}
  A minimal Boolean $\pm$\hyp{}preclone $F$ has minimal signum
  $(\mathord{+}, \allowbreak \mathord{+}, \allowbreak \mathord{+}, \allowbreak \mathord{-})$
  if and only if $F$ is
  generated by one of the following two $\pm$\hyp{}operations with minimal
  signum $(\mathord{+},\mathord{+},\mathord{+},\mathord{-})$,
  $f_{\mathrm{maj}}$ \textup{(}majority operation on variables with
  sign $+$\textup{)} 
  or $f_{\oplus}$ \textup{(}minority operation on variables with sign
  $+$\textup{)}:  
\begin{align*}
        f_{\mathrm{maj}}(x,y,z,u)&:=(x\land y)\lor (y\land z)\lor (z\land x),
         & f_{\oplus}(x,y,z,u):=x\oplus y \oplus z.
      \end{align*}
\end{lemma}

\begin{proof}
  From Lemma \ref{lem:nearprojection} we can conclude that
  if $F$ is a minimal Boolean $\pm$\hyp{}preclone with minimal signum
  $(\mathord{+},\mathord{+},\mathord{+},\mathord{-})$,
  then $F$ contains a $\pm$\hyp{}near\hyp{}projection $f$ of signum
  $(\mathord{+},\mathord{+},\mathord{+},\mathord{-})$ whose last
  argument (with sign $-$) is fictitious. Because $F$ is generated by
  $f$, it follows that each argument with sign $-$ of each
  $\pm$\hyp{}operation in $F$ must be fictitious and $\mathring{F}$
  essentially coincides with $\Sg{\mathring{F}}$. Thus
  $\Sg{\mathring{F}}$ must be a minimal Boolean clone generated by the
  essentially ternary operation $\mathring{f}$.
 Conversely, if a ternary operation $g$ generates a minimal clone,
 then the $4$\hyp{}ary $\pm$\hyp{}operation $f$ with $\sgn(f)=(\mathord{+},\mathord{+},\mathord{+},\mathord{-})$ given
 by $f(x,y,z,u):=g(x,y,z)$ generates a minimal
 $\pm$\hyp{}preclone. It is known that there are exactly two minimal
 Boolean clones
 generated by a ternary (but not by a binary) operation $g$, which leads to
 the two possible generators $f_{\mathrm{maj}}$ and $f_{\oplus}$ for $F$.
\end{proof}

\begin{notation}
  In what follows we often have to deal with binary Boolean functions
  $\alpha \colon \{0,1\}^{2}\to \{0,1\}$ which we shall enumerate via the
  notation $\alpha_{i_{1}i_{2}i_{3}i_{4}}$ with $i_{1}=\alpha(0,0)$,
  $i_{2}=\alpha(0,1)$, $i_{3}=\alpha(1,0)$, $i_{4}=\alpha(1,1)$.
  For brevity, we may also refer to $\alpha$ simply by the $0$--$1$ sequence
  $\alpha = i_{1}i_{2}i_{3}i_{4} = \alpha(0,0) \alpha(0,1) \alpha(1,0) \alpha(1,1)$,
  especially when $\alpha$ appears as an index in notations such as
  $h_{\alpha}$ and $g_{\alpha}$ in the next lemmas.
\end{notation}

\begin{lemma}\label{min-4aryB}
\leavevmode
  \begin{enumerate}[label=\textup{(\roman*)}]
  \item\label{4aryBi}
  Let $F$ be a minimal Boolean $\pm$\hyp{}preclone with minimal signum
  $(\mathord{+},\mathord{+},\mathord{-},\mathord{-})$.
  Then $F$ is generated by a $\pm$\hyp{}near\hyp{}projection $h$
  with $\sgn(h)=(\mathord{+},\mathord{+},\mathord{-},\mathord{-})$ satisfying $h(x,x,y,z)=x=h(x,y,z,z)$. More
  explicitly, $h$ can be expressed in the form
  \begin{align*}
    h_{\alpha}(x,y,z,u)=
    \begin{cases}
      x &\text{if } x=y \text{ or }z=u\\
      \alpha(x,z) &\text{otherwise, i.e., } x\neq y \text{ and }z\neq
      u,
    \end{cases}
  \end{align*}
  where $\alpha \colon \{0,1\}^{2}\to\{0,1\}$ is a Boolean function except
  the first projection $\alpha_{0011}(x,z)=x$.

  \item\label{4aryBii}
  We have $h_{\alpha'}(x,y,z,u)=h_{\alpha}(x,y,u,z)$ for
  $\alpha'(x,z):=\alpha(x,\lnot z)$ \textup{(}consequently
  $\pmSg{h_{\alpha'}}=\pmSg{h_{\alpha}}$\textup{)}. Moreover, for dual
  operations we have $(h_{\alpha})^{\partial} = h_{\alpha^{\partial}}$.

  \item\label{4aryBiii}
  $\pmSg{h_{0000}} = \pmSg{h_{0001}} = \pmSg{h_{0010}}$.

  \item\label{4aryBiv}
  Each of the Boolean $\pm$\hyp{}operations $h_{0000}$, $h_{1111}$ and $h_{1100}$
  has minimal signum $(+,+,-,-)$.

  \end{enumerate}
\end{lemma}

\begin{proof}
  \ref{4aryBi}: Since $h$ is a $\pm$\hyp{}near\hyp{}projection we have
  $h(x,x,z,u) = x$ and $h(x,y,z,z) \allowbreak  = x$ or $h(x,y,z,z) = y$, in the latter
  case we may exchange $x$ and $y$ to get a $\pm$\hyp{}near\hyp{}projection satisfying
  $h(x,x,z,u)=x=h(x,y,z,z)$.

  Note that $x \neq y \iff y = \lnot x$. Thus, with
  $\alpha(x,z) := h(x,\lnot x,z,\lnot z)$, we get $h=h_{\alpha}$ as
  defined above.

  \ref{4aryBii}:
  $h_{\alpha'}(x,\lnot x,z,\lnot z)=\alpha'(x,z)=\alpha(x,\lnot
  z)=h_{\alpha}(x,\lnot x,\lnot z,z)$ gives
  $h_{\alpha'}(x,y,z,u) \allowbreak = h_{\alpha}(x,y,u,z)$.
  The statement $(h_{\alpha})^{\partial} = h_{\alpha^{\partial}}$ is
  clear.

\ref{4aryBiii}: In fact,\begin{align*}
   h_{0000}(x,y,z,u) &= h_{0001}(h_{0001}(x,y,z,u),y,u,z), \\
   h_{0001}(x,y,z,u) &= h_{0010}(x,y,u,z) \qquad \text{(since $(0001)'=(0010)$)}, \\
   h_{0010}(x,y,z,u) &= h_{0000}(x,y,h_{0000}(u,z,x,y),z).
 \end{align*}

\ref{4aryBiv}:
  We are going to show that $h_{0000}$ and $h_{1100}$ have minimal signum. Then
  also $h_{1111}$ has this minimal signum since
  $h_{1111}=h_{0000^{\partial}}=h_{0000}^{\partial}$ (cf.\ Lemma~\ref{lem:partial}).

Because each $\pmSg{h_{\alpha}}$ must contain a minimal $\pm$\hyp{}preclone
generated by a quaternary $\pm$\hyp{}near projection
(see Lemmas \ref{lem:nearprojection} and
\ref{lem:Snearprojection1}), we simply check that neither $h_{0000}$ nor
$h_{1100}$ generates neither $f_{\mathrm maj}$ nor $f_{\oplus}$
nor any $h_{\alpha}$ which generates another
$\pm$\hyp{}preclone. In order to show that a $\pm$\hyp{}operation $h$ does not
generate some $f$
we exhibit a $\pm$\hyp{}relation $S$\hyp{}preserved by $h$ but not by $f$.

Concerning $h_{0000}$: Let $\rho:=(\rho_{+},\rho_{-})$ with
$\rho_{+}:=\left(
\begin{smallmatrix}
  0&0&1\\ 0&1&0
\end{smallmatrix}\right)$
and $\rho_{-}:=\left(
\begin{smallmatrix}
  0&0&1&1\\ 0&1&0&1
\end{smallmatrix}\right)$.
It is straightforward to check that
$h_{0000}\Spreserves\rho$. Therefore each
$h_{\alpha}\in\pmSg{h_{0000}}$ must $S$\hyp{}preserve $\rho$, too.
If $\alpha(0,0)=1$ then $h_{\alpha}(\pair 01,\pair 10,\pair
00,\pair 10)=\pair 11\notin\rho_{+}$, which shows that $h_{\alpha}$ does not
$S$\hyp{}preserve $\rho$. Analogously, $\alpha(0,1)=1$ implies
$h_{\alpha}(\pair 01,\pair 10,\pair 10,\pair 00)=\pair 11$ showing
that $h_{\alpha}$ does not $S$\hyp{}preserve $\rho$. There remain the
$\pm$\hyp{}near\hyp{}projections $h_{\alpha}$ with
$\alpha\in\{\alpha_{0011},\alpha_{0001},\alpha_{0010},\alpha_{0000}\}$.

We dismiss $h_{0011}$, as it is a trivial projection.
For the remaining operations we have
 $\pmSg{h_{0000}}=\pmSg{h_{0001}}=\pmSg{h_{0010}}$ as shown above in \ref{4aryBiii}.
Furthermore, let the ternary $\pm$\hyp{}relation $\bar\sigma$ be given by $\bar\sigma_{+}:=\{0,1\}^{3}\setminus \{(1,1,1)\}$ and
$\bar\sigma_{-}:=\{(0,0,0)\}$. Then $h_{0000}\Spreserves\bar\sigma$
but $f_{\mathrm{maj}}$ and $f_{\oplus}$ do not preserve $\bar\sigma$
since
$f_{\mathrm{maj}}(\triple 011,\triple 101,\triple 110,\triple 000) \allowbreak =
\triple 111 \notin\bar\sigma_{+}$ and $f_{\oplus}(\triple 100,\triple
010,\triple 001,\triple 000)=\triple 111\notin\bar\sigma_{+}$.

Concerning $h_{1100}$, note that $h_{1100}$ is self-dual since
$1100^{\partial}=1100$. Thus each
$4$\hyp{}ary near\hyp{}projection $h_{\alpha}$ in $\pmSg{h_{1100}}$ must be
self-dual, too. Except the trivial projection $h_{0011}$, there are only
two such near\hyp{}projections, namely $h_{0101}$ and $h_{1010}$. Moreover,
they generate the same $\pm$\hyp{}preclone because
$\alpha_{1010}=\alpha_{0101}'$. We show
$\pmSg{h_{1100}}\subseteq \pmSg{h_{0101}}$, which implies that $\pmSg{h_{1100}}$
cannot properly contain a $\pm$\hyp{}preclone generated by some other
$h_{\alpha}$. In fact, it is easy to verify that
$h_{1100}(x,y,z,u)=h_{0101}(h_{0101}(x,y,z,u),y,u,z)$.
Moreover, $h_{0101}\Spreserves\bar\sigma$ but $f_{\mathrm maj}$ and
$f_{\oplus}$ do not $S$\hyp{}preserve $\bar\sigma$ (as already shown above)
and we are done.
\end{proof}

\begin{lemma}\label{min-3aryA}
\leavevmode
  \begin{enumerate}[label=\textup{(\roman*)}]
  \item\label{3aryAi}
  Let $F$ be a minimal Boolean $\pm$\hyp{}preclone with minimal signum
  $(\mathord{+},\mathord{+},\mathord{-})$. Then $F$ is generated by a $\pm$\hyp{}near\hyp{}projection $g$
  with $\sgn(g)=(\mathord{+},\mathord{+},\mathord{-})$ satisfying $g(x,x,y)=x$. Thus $g$ can be
  expressed in the form $g=g_{\alpha}$ with
  \begin{align*}
    g_{\alpha}(x,y,z)=
    \begin{cases}
      x &\text{if } x=y \\
      \alpha(x,z) &\text{otherwise},
    \end{cases}
  \end{align*}
  where $\alpha \colon \{0,1\}^{2}\to\{0,1\}$ is a Boolean function except
  $\alpha_{0011}(x,z)=x$ or $\alpha_{1100}(x,z)=\lnot x$ \textup{(}then
  $g_{\alpha}$ is the first or second projection, respectively\textup{)}.

  \item\label{3aryAii}
  We have $g_{\alpha^{*}}(x,y,z)=g_{\alpha}(y,x,z)$ for
  $\alpha^{*}(x,z):=\alpha(\lnot x,z)$ \textup{(}consequently
  $\pmSg{g_{\alpha^{*}}}=\pmSg{g_{\alpha}}$\textup{)}. Moreover we have
  $(g_{\alpha})^{\partial}=g_{\alpha^{\partial}}$.

  \item\label{3aryAiii}
  $g_{\alpha}$ has minimal signum $(\mathord{+},\mathord{+},\mathord{-})$ for each nontrivial
  $\alpha$ except $\alpha_{0110}$ and $\alpha_{1001}$.

  \item\label{3aryAiv}
  For $\alpha \neq \beta$,
  neither of
  $\pmSg{g_{\alpha}}$ and $\pmSg{g_{\beta}}$ is
  contained in the other unless $\beta=\alpha^{*}$.

\end{enumerate}
\end{lemma}

\begin{proof} \ref{3aryAi}: For given $g$ satisfying $g(x,x,y)=x$, define
  $\alpha(x,z):=g(x,\lnot x,z)$. Then $g=g_{\alpha}$. Conversely, each
  $g_{\alpha}$ satisfies $g_{\alpha}(x,x,y)=x$ by definition. 

\ref{3aryAii}: 
$g_{\alpha^{*}}(x,\lnot x,z)=\alpha^{*}(x,z)=\alpha(\lnot
x,z)=g_{\alpha}(\lnot x,\lnot \lnot x,z)=g_{\alpha}(\lnot x,x,z)$ shows that
$g_{\alpha^{*}}(x,y,z)=g_{\alpha}(y,x,z)$.
Moreover,
$g_{\alpha}^{\partial}(x,\lnot x,z)
=\lnot(g_{\alpha}(\lnot x,x,\lnot z))
=\lnot(\alpha(\lnot x,\lnot z))
=\alpha^{\partial}(x,z)=g_{\alpha^{\partial}}(x,\lnot x,z)$
implies
$(g_{\alpha})^{\partial}=g_{\alpha^{\partial}}$.

\ref{3aryAiii}: The connections between the mappings $\alpha$,
  $\alpha^{*}$ and $\alpha^{\partial}$ are the following:
  \begin{align*}
\alpha_{0000}&=\alpha_{0000}^{*}=\alpha_{1111}^{\partial},&
                                              \alpha_{0001}
          &=\alpha_{0100}^{*}=\alpha_{0111}^{\partial}
            =(\alpha_{1101}^{\partial})^{*},\\
\alpha_{0110}&=\alpha_{0110}^{\partial}=\alpha_{1001}^{*},
        & \alpha_{0010}&=\alpha_{1000}^{*}
 =\alpha_{1011}^{\partial}=(\alpha_{1110}^{\partial})^{*},\\
    (\text{note }
\alpha_{0101}&=\alpha_{0101}^{*}=\alpha_{0101}^{\partial},&
\alpha_{1010}&=\alpha_{1010}^{*}=\alpha_{1010}^{\partial} \text{).}
  \end{align*}
We are going to show that $g_{\alpha}$ has minimal signum
  $(\mathord{+},\mathord{+},\mathord{-})$ for
  $\alpha \in \{\alpha_{0000}, \allowbreak \alpha_{0001}, \allowbreak \alpha_{0010}, \allowbreak \alpha_{0101}, \allowbreak \alpha_{1010}\}$.
Because
$\pmSg{g_{\alpha}}=\pmSg{g_{\alpha^{*}}}$ is minimal if and only if
$\pmSg{g_{\alpha}^{\partial}}$ is minimal by Lemma~\ref{lem:partial},
this shows minimality also for all other $g_{\alpha}$ except $\alpha_{0110}$
and $\alpha_{1001}$ which will be treated in the proof of Theorem
\ref{min-theorem}\ref{min-thmC} by showing that
$\pmSg{g_{0110}}=\pmSg{g_{1001}}$ is not a
minimal $\pm$\hyp{}clone.

For these five $\alpha$ we prove first that they do not generate a
$4$\hyp{}ary near\hyp{}projection.

Let $\sigma:=(\sigma_{+},\sigma_{-})$ with $\sigma_{+}:=
\left(\begin{smallmatrix}
  0&0&1\\0&1&1
\end{smallmatrix}\right)
$
and $\sigma_{-}:=\left(
\begin{smallmatrix}
  0&0&1&1\\ 0&1&0&1
\end{smallmatrix}\right)$.
Then it is straightforward to prove that
$g_{\alpha}\Spreserves\sigma$ for each of the five functions $\alpha$ above.
However, no nontrivial $h_{\beta}$ (cf.\ Lemma~\ref{min-4aryB}) $S$\hyp{}preserves
$\sigma$; in fact, for each $\beta$ except $\beta=0011$ (which gives
trivial $h_{\beta}$) there exist
$c\neq d$ in $\{0,1\}$ such that $h_{\beta}(0,1,c,d)=1$ or
$h_{\beta}(1,0,c,d)=0$. Then
$h_{\beta}(\pair 00,\pair 11,\pair c1,\pair d1)=\pair
10\notin\sigma_{+}$ or $h_{\beta}(\pair 11,\pair 00,\pair 0c,\pair 0d)=\pair
10\notin\sigma_{+}$, respectively, shows that $h_{\beta}$ does not
$S$\hyp{}preserve $\sigma$.

Furthermore, neither does $f_{\oplus}$ preserve $\sigma$, since
$f_{\oplus}(\pair 00,\pair 01,\pair 11,\pair 00)=\pair 10\notin
\sigma_{+}$.

Finally, $g_{\alpha}$ preserves the following ternary $\bar\sigma$ (as
in the proof of Lemma~\ref{min-4aryB}\ref{4aryBiv}) or $\hat\sigma$ (we omit
the straightforward verification) but
$f_{\mathrm{maj}}$ does not:

For $\alpha\in\{\alpha_{0000},\alpha_{0001},\alpha_{0101}\}$, take
$\bar\sigma_{+}:=\{0,1\}^{3}\setminus \{(1,1,1)\}$ and
$\bar\sigma_{-}:=\{(0,0,0)\}$.

For $\alpha\in\{\alpha_{0010},\alpha_{1010}\}$, take
$\hat\sigma_{+}:=\bar\sigma_{+}$
and
$\hat\sigma_{-}:=\{(1,1,1)\}$.

Then
$f_{\mathrm{maj}}(\triple 011,\triple 101,\triple 110,\triple
000)=f_{\mathrm{maj}}(\triple 011,\triple 101,\triple 110,\triple
111)=\triple 111\notin\bar\sigma_{+}=\hat\sigma_{+}$.

Consequently, for each of the above $g_{\alpha}$, the $\pm$\hyp{}preclone
$\pmSg{g_{\alpha}}$
cannot contain a $4$\hyp{}ary near\hyp{}projection. Therefore each minimal $\pm$\hyp{}preclone contained in
$\pmSg{g_{\alpha}}$, must contain some near\hyp{}projection $g_{\beta}$.

\ref{3aryAiv}:

For
$\alpha=\alpha_{abcd}$ let
$\rho^{\alpha}:=(\left[
  \begin{smallmatrix}
    0&1&a&c\\
    0&1&b&d\\
    1&0&c&a\\
    1&0&d&b
  \end{smallmatrix}
\right],
\left[\begin{smallmatrix}
  0\\1\\0\\1
\end{smallmatrix}\right])$.
Then $\alpha^{*}=\alpha_{cdab}$ and $\rho^{\alpha}=\rho^{\alpha^{*}}$.
Furthermore, $g_{\alpha}\Spreserves\rho^{\alpha}$ (easy exercise) but
obviously $g_{\beta}$ does
not preserve $\rho^{\alpha}$ for $\beta\notin\{\alpha,\alpha^{*}\}$.
This shows that $\pmSg{g_{\beta}}$ is not contained in $\pmSg{g_{\alpha}}$
and we are done.
\end{proof}

\begin{lemma}\label{min-3aryB}
There are no $\pm$\hyp{}operations with minimal signum $(\mathord{+},\mathord{-},\mathord{-})$.
\end{lemma}

\begin{proof}
Assume that $f$ generates a minimal $\pm$\hyp{}preclone
  $F$ such that $\sgn(f)=(\mathord{+},\mathord{-},\mathord{-})$ is the minimal signum of $F$. Then
  $f(x,y,y)=x$ since otherwise $g(x,y):=f(x,y,y)$ would be nontrivial
  with $\sgn(g)=(\mathord{+},\mathord{-})$ and $\sgn(f)$ would not be the minimal signum of
  $F$. Thus there exists a binary Boolean function $\beta$ (just take
  $\beta(x,y):=f(x,y,\lnot y)$) which is not the first projection
  $\alpha_{0011}$, such that
  \begin{align*}
    f(x,y,z)=f^{\beta}(x,y,z):=
    \begin{cases}
      x &\text{if } y=z,\\
      \beta(x,y) &\text{otherwise.}
    \end{cases}
  \end{align*}
For example,
  \begin{align*}
    f^{1100}(x,y,z)=x\oplus y\oplus z=
    \begin{cases}
      x &\text{if }{y=z},\\
      \lnot x &\text{otherwise.}
    \end{cases}
  \end{align*}

At first we exclude $f^{1100}$: We have $f^{1100}(x,f^{1100}(u,y,z),u)=x\oplus(u\oplus y\oplus
z)\oplus u=x\oplus y\oplus z=f_{\oplus}(x,y,z,u)$. The $\pm$\hyp{}operation
$f_{\oplus}$ has minimal signum $(\mathord{+},\mathord{+},\mathord{+},\mathord{-})$ by Lemma
\ref{min-4aryA}. Therefore $f^{1100}$ cannot generate a minimal $\pm$\hyp{}preclone.

  Now we show that for all remaining $f^{\beta}$ (with
  $f^{\beta}\notin\{f^{0011},f^{1100}\}$) there exists a nontrivial
  $\pm$\hyp{}near\hyp{}projection in $\pmSg{f^{\beta}}$ of signum $(\mathord{+},\mathord{+},\mathord{-})$
  (with $2>1$ signs `$+$'), which shows that $(\mathord{+},\mathord{-},\mathord{-})$ cannot be the
  minimal signum of $\pmSg{f^{\beta}}$ (then one must get a
  $\pm$\hyp{}preclone in the list in Theorem~\ref{min-theorem}\ref{min-thmC}).

  Let $g^{\beta}(x,y,z):=f^{\beta}(x,z,f^{\beta}(z,x,y))$ and
  $\beta'(x,y):=\beta(x,\lnot y)$. Then
  $\sgn(g^{\beta})=(+,+,-)$. Moreover,
  $f^{\beta'}(x,y,z)=f^{\beta}(x,z,y)$ since, for $y\neq z$, we have
  $f^{\beta'}(x,y,z) \allowbreak =f^{\beta'}(x,\lnot z,z)=\beta'(x,\lnot
  z)=\beta(x,z)=g^{\beta}(x,z,\lnot z)=f^{\beta}(x,z,y)$, in
  particular $\pmSg{g^{\beta'}}=\pmSg{g^{\beta}}$.

  Note that
  $g^{\beta}$ is a $\pm$\hyp{}near\hyp{}projection since
  $g^{\beta}(x,x,z)=f^{\beta}(x,z,f^{\beta}(z,x,x))=f^{\beta}(x,z,z)=x$. We
  are going to show that $g^{\beta}$ is nontrivial, i.e., neither the
  first nor second projection.

  We compute
$g^{\beta}(0,1,0)=f^{\beta}(0,0,f^{\beta}(0,0,1))=f^{\beta}(0,0,\beta(0,0))=\beta(0,0)$
  (note that $f^{\beta}(0,0,0)=0$ and $f^{\beta}(0,0,1)=\beta(0,0)$), analogously
  $g^{\beta}(1,0,1)=\beta(1,1)$. Thus $g^{\beta}$ is not the first projection
  if $\beta(0,0)=1$ or $\beta(1,1)=0$ (since then $g^{\beta}(0,1,0)=1$
  or $g^{\beta}(1,0,1)=0$), and not the second projection if
  $\beta(0,0)=0$ or $\beta(1,1)=1$. Therefore each of the following
  $f^{\beta}$ does not have minimal signum $(\mathord{+},\mathord{-},\mathord{-})$:
  $f^{1ij1}, f^{0ij0}$ as well as the corresponding $f^{\beta'}$,
  i.e., $f^{i11j},f^{i00j}$ ($i,j\in\{0,1\}$).

  The only missing $f^{\beta}$ are $f^{ijj'i'}$ with $i\neq i'$,
  $j\neq j'$, i.e., $f^{0011},f^{1100}$ (already excluded above) and
  $f^{\beta_{0}}:=f^{0101}$ and $ f^{\beta_{0}'}=f^{1010}$. Thus it
  remains to exclude
  \begin{align*}
    f^{\beta_{0}}(x,y,z)=
    \begin{cases}
      x &\text{if } y=z,\\
      y &\text{otherwise.}
    \end{cases}
  \end{align*}
  For this, we observe that $g_{0}(x,y,z):=f^{\beta_{0}}(x,f^{\beta_{0}}(z,y,x),z)$ is a nontrivial
  $\pm$\hyp{}near\hyp{}projection with signum $(\mathord{+},\mathord{+},\mathord{-})$. In fact, we have
  $g_{0}(x,x,z)=x$, $g_{0}(0,1,0)=f^{\beta_{0}}(0,f^{\beta_{0}}(0,1,0),0)=f^{\beta_{0}}(0,1,0)=1$ and
$g_{0}(1,0,0)=f^{\beta_{0}}(1,f^{\beta_{0}}(0,0,1),0)=f^{\beta_{0}}(1,0,0)=1$,
  showing that $g_{0}$ is
  not a trivial $\pm$\hyp{}projection.
\end{proof}

\begin{lemma}\label{min-2ary}
  Each of the following five binary Boolean functions considered as a
  $\pm$\hyp{}operation with signum $(\mathord{+},\mathord{-})$ generates a minimal
  $\pm$\hyp{}preclone with minimal signum $(\mathord{+},\mathord{-})$:
  $\alpha_{0000}=\alpha_{1111}^{\partial}$,
  $\alpha_{0001}=\alpha_{0111}^{\partial}$, $\alpha_{0101}$,
  $\alpha_{1010}$, $\alpha_{1100}$.
  Moreover, each of these five $\pm$\hyp{}operations is the unique $\pm$\hyp{}operation of signum $(\mathord{+},\mathord{-})$ in the $\pm$\hyp{}preclone it generates.
\end{lemma}

\begin{proof} For each $\pm$\hyp{}preclone $F=\pmSg{\alpha}$ (with
  $\alpha$ as listed in the lemma), it is easy to see that for each
  nontrivial term operation $f \in F$ identifying all variables of the
  same sign shall give back the generating $\pm$\hyp{}operation $\alpha$.
  Therefore, $F$ is minimal.

To see this, we just give a few hints and leave the details to the reader.
For $\alpha_{0000}(x,y)=0$, all functions
in $F$ are constant $0$. For $\alpha_{0101}(x,y)=y$,
any term operation $f(x_{1},\dots,x_{n})$ equals a projection corresponding to the rightmost variable
in the term. For $\alpha_{1010}(x,y)=\lnot y$ and a term operation
$f(x_{1},\dots,x_{n})$, assume $x_{n}$ is the rightmost variable in the term. If
$x_{n}$ has sign $+$, then the term operation equals $x_{n}$ and is a
trivial projection, otherwise it equals $\lnot x_{n}$. For
$\alpha_{1100}(x,y)=\lnot x$, each term function
$f(x_{1},\dots,x_{n})$ equals the leftmost variable (say $x_{1}$) or
its negation $\lnot x_{1}$ depending on the depth of the variable in
the term.
Finally, for $\alpha_{0001}(x,y)=x\land y$, each term function
$f(x_{1},\dots,x_{n})$ is the conjunction of some variables (and there
are always variables of sign $+$ as well as of sign $-$ among them).
\end{proof}

Now we are ready to characterize all minimal Boolean $\pm$\hyp{}preclones.

  \begin{theorem}[\textbf{minimal Boolean $\pm$\hyp{}preclones}]\label{min-theorem}
    There are twenty minimal Boolean $\pm$\hyp{}preclones. They are of
    the form $F=\pmSg{f}$ where $f$ is one of the $\pm$\hyp{}operations in
    the following list:
    \begin{enumerate}[label=\textup{(\Alph*)}]

    \item\label{min-thmA} \emph{minimal signum $\lambda=(\mathord{+},\mathord{+},\mathord{+},\mathord{-})$:}
      There are two minimal $\pm$\hyp{}preclones with minimal signum
      $\lambda$, generated by $f_{\mathrm{maj}}$ or $f_{\oplus}$
      \textup{(}$\sgn(f)=\lambda$, the last argument is fictitious\textup{)} where
      \textup{(}cf.\ Lemma~\ref{min-4aryA}\textup{)}
      \begin{align*}
        f_{\mathrm{maj}}(x,y,z,u)&:=(x\land y)\lor (y\land z)\lor
                        (z\land x), && \text{\textup{(majority function}} \\ &&& \text{\textup{on
                                   $x,y,z$)} }\\
 f_{\oplus}(x,y,z,u)&:=x\oplus y\oplus z.
      \end{align*}

    \item\label{min-thmB} \emph{minimal signum $\lambda=(\mathord{+},\mathord{+},\mathord{-},\mathord{-})$:}
      There are three minimal $\pm$\hyp{}preclones with minimal signum
      $\lambda$, generated by $f\in\{h_{0000},h_{1111},h_{1100}\}$
      \textup{(}$\sgn(f)=\lambda$\textup{)} where \textup{(}cf.\ Lemma~\ref{min-4aryB}\textup{)}
      \begin{align*}
        h_{0000}(x,y,z,u)
            &= x\land(y\lor(z\land u)\lor(\lnot z\land\lnot u))
         \\ &=
                           \begin{cases}
                             x &\text{if $x=y$ or $z=u$,} \\
                             0 & \text{otherwise,}
                           \end{cases}\\
        h_{1111}(x,y,z,u)
            &= x\lor(y\land(z\lor u)\land(\lnot z\lor\lnot u))
         \\ &=
                           \begin{cases}
                             x &\text{if $x=y$ or $z=u$,} \\
                             1 & \text{otherwise,}
                           \end{cases}\\
        h_{1100}(x,y,z,u)&=
                          ( x \land y ) \lor ( x \land ( z \oplus \lnot u ) ) \lor ( y \land ( z \oplus u ) )
                         \\
                         &=
                           \begin{cases}
                             x&\text{if $x=y$ or $z=u$,} \\
                             \lnot x\;(=y) & \text{otherwise.}
                           \end{cases} 
      \end{align*}

    \item\label{min-thmC} \emph{minimal signum $\lambda=(\mathord{+},\mathord{+},\mathord{-})$:}
      There are eight minimal $\pm$\hyp{}preclones with minimal signum
      $\lambda$, generated by a $\pm$\hyp{}near\hyp{}projection $f$
      \textup{(}$\sgn(f)=\lambda$\textup{)} where $f$ is one of the following
      $\pm$\hyp{}operations \textup{(}cf.\ Lemma~\ref{min-3aryA}\textup{)}:
\begin{align*}
  g_{0000}(x,y,z)&=x\land y=
                   \begin{cases}
                     x &\text{if $x=y$} \\
                     0 &\text{otherwise,}
                   \end{cases}\\
  g_{1111}(x,y,z)&=x\lor y=
                   \begin{cases}
                     x &\text{if $x=y$} \\
                     1 &\text{otherwise,}
                   \end{cases}
\end{align*}\begin{align*}
g_{0101}(x,y,z)&=(x\land y)\lor (y\land z)\lor (z\land x)=
                   \begin{cases}
                     x &\text{if $x=y$} \\
                     z &\text{otherwise,}
                   \end{cases}\\
  g_{1010}(x,y,z)&=(x\land y)\lor (y\land\lnot z)\lor (\lnot z\land
                   x)=
                   \begin{cases}
                     x &\text{if $x=y$} \\
                     \lnot z &\text{otherwise,}
                   \end{cases}
\end{align*}\begin{align*}
  g_{0111}(x,y,z)&=x\lor (y\land z)=
                   \begin{cases}
                     x &\text{if $x=y$} \\
                     x\lor z &\text{otherwise,}
                   \end{cases}\\
  g_{0001}(x,y,z)&=x\land (y\lor z)=
                   \begin{cases}
                     x &\text{if $x=y$} \\
                     x\land z &\text{otherwise,}
                   \end{cases}
\end{align*}\begin{align*}
  g_{1011}(x,y,z)&=x\lor (y\land\neg z)=
                   \begin{cases}
                     x &\text{if $x=y$} \\
                     x \lor \lnot z &\text{otherwise,}
                   \end{cases}\\
  g_{0010}(x,y,z)&=x\land (y\lor\neg z)=
                   \begin{cases}
                     x &\text{if $x=y$} \\
                     x\land \lnot z &\text{otherwise.}
                   \end{cases}
\end{align*}

\item\label{min-thmD} \emph{minimal signum $\lambda=(\mathord{+},\mathord{-})$:}
  There are
  seven minimal $\pm$\hyp{}preclones with minimal signum $\lambda$,
  generated by one of the following $\pm$\hyp{}operations \textup{(}with signum
  $\lambda$, cf.\ Lemma~\ref{min-2ary}\textup{)}:
  \begin{align*}
    \alpha_{0000}(x,y)&=0,& \alpha_{1111}(x,y)&=1,  \\
    \alpha_{0001}(x,y)&=x\land y, & \alpha_{0111}(x,y)&=x\lor y,\\
    \alpha_{0101}(x,y)&=y,& \alpha_{1010}(x,y)&=\lnot y,
                      &\alpha_{1100}(x,y)=\lnot x.
  \end{align*}
\end{enumerate}

\end{theorem}

\begin{proof} By Lemmas \ref{lem:nearprojection} and
  \ref{min-3aryB} there exist exactly the 4 minimal signa
  treated in \ref{min-thmA}--\ref{min-thmD}. Lemmas
  \ref{min-4aryA}, \ref{min-4aryB}, \ref{min-3aryA} and \ref{min-2ary} show that all listed $\pm$\hyp{}operations
  generate a mininal $\pm$\hyp{}preclone and have minimal signum. Thus it
  remains to prove that all other possible candidates for generators
  of a minimal $\pm$\hyp{}preclone do not generate a minimal preclone or
  generate a $\pm$\hyp{}operation that is already included in the list.

\ref{min-thmA}: This was already done in Lemma~\ref{min-4aryA}.

  \ref{min-thmB}:
  We show that each nontrivial $\pm$\hyp{}operation $h_{\alpha}$ with
  signum $(\mathord{+},\mathord{+},\mathord{-},\mathord{-})$ generates a $\pm$\hyp{}operation $f_{\alpha}$ which is
  one of the $\pm$\hyp{}operations $h_{0000}$, $h_{1111}$, $h_{1100}$ with
  minimal signum from Lemma~\ref{min-4aryB}. According
  to Lemma~\ref{min-4aryB} only one of the Boolean functions $\alpha$ and
  $\alpha'$ must be checked.

  By Lemma~\ref{min-4aryB}\ref{4aryBiii},
  we know that $h_{0000}$, $h_{0001}$, and $h_{0010}$ generate the same
  (minimal) $\pm$\hyp{}preclone. Thus also the dual
  $\pm$\hyp{}operations $h_{1111}$, $h_{0111}$, and $h_{1011}$ generate the
  same minimal $\pm$\hyp{}preclone.

It remains to check four cases
$\alpha\in\{\alpha_{1000},\alpha_{1101},\alpha_{0110},\alpha_{0101}\}$
  (note that then
$\alpha'\in\{\alpha_{1110},\alpha_{0100},\alpha_{1001},\alpha_{1010}\}$
and therefore all candidates for $\alpha$ are covered). For these four cases let
$f_{\alpha}(x,y,z,u):=h_{\alpha}(h_{\alpha}(x,y,z,u),y,u,z)$. Then it
is straightforward to check $f_{\alpha}(x,y,z,u)=h_{1100}(x,y,z,u)$,
the latter generating a minimal $\pm$\hyp{}preclone, and we are done.

  \ref{min-thmC}: At first we show that
  \begin{align*}
  g_{0110}(x,y,z)&=(x\land\lnot z)\lor (y\land z)=
                   \begin{cases}
                     x &\text{if $x=y$} \\
                     (x\land \lnot z)\lor (\lnot x \land
                     z)&\text{otherwise}
                   \end{cases}
  \end{align*}
(and therefore also $g_{1001}$ because
$\alpha_{1001}=\alpha_{0110}^{*}$) does \textbf{not} generate a
minimal $\pm$\hyp{}preclone. In fact, we have
$g_{0110}(g_{0110}(x,y,u),g_{0110}(y,x,u),z)=h_{1100}(x,y,z,u)$

(cf.\ \ref{min-thmB}).

By Lemma~\ref{min-3aryA},
  all remaining $\pmSg{g_{\alpha}}$ are minimal $\pm$\hyp{}preclones and
  $\pmSg{g_{\alpha}}=\pmSg{g_{\beta}}$ if and only if
  $\beta\in\{\alpha,\alpha^{*}\}$. With the help of the list of all $\alpha$ and
  $\alpha^{*}$ given in the proof of Lemma \ref{min-3aryA}, we see that there remain
  exactly the eight
  different $\pm$\hyp{}preclones listed in \ref{min-thmC} that are
  generated by a $\pm$\hyp{}operation with
  minimal signum $(\mathord{+},\mathord{+},\mathord{-})$.

  \ref{min-thmD}: According to Lemma~\ref{min-2ary} all seven listed
  $\pm$\hyp{}operations generate a minimal $\pm$\hyp{}preclone (taking into
  account also the dual operations) and all these $\pm$\hyp{}preclones are
  different (since $\alpha$ is the only nontrivial operation with
  signum $(\mathord{+},\mathord{-})$ in $\pmSg{\alpha}$, see Lemma
  \ref{min-2ary}). Now we show that there are no further minimal $\pm$\hyp{}preclones with minimal signum $(\mathord{+},\mathord{-})$.

  In fact, we show for each of the remaining eight nontrivial Boolean
  $\pm$\hyp{}operations $\alpha \colon \{0,1\}^{2} \to \{0,1\}$ considered as a
  $\pm$\hyp{}operation with signum $(\mathord{+},\mathord{-})$ that $\pmSg{\alpha}$ contains an
  already listed $\pm$\hyp{}operation $f$ with minimal $\pm$\hyp{}preclone
  $\pmSg{f}$. In the following we indicate $f$ and the term that
  produces $f$ (it is an easy exercise to check the correctness):

  $\alpha$ equals $\alpha_{0110}(x,y)=x\oplus y$ or
  $\alpha_{1001}(x,y)=1\oplus x\oplus y$:
  \newline
  Then
  $f(x,y):=\alpha(x,\alpha(y,x))= y=\alpha_{0101}(x,y)$.

  $\alpha$ equals $\alpha_{1101}(x,y)=\lnot x\lor y$:
  \newline
  Then
  $f(x,y)=\alpha(x,\alpha(y,x))=1=\alpha_{1111}(x,y)$.

  $\alpha$ equals $\alpha_{0100}(x,y)=\lnot x\land y$:
  \newline
  Then
  $f(x,y)=\alpha(x,\alpha(y,x))=0=\alpha_{0000}(x,y)$.

  $\alpha$ equals $\alpha_{1011}(x,y)=x\lor \lnot y$:
  \newline
  Then
  $f(x,y,z)=\alpha(x,\alpha(z,y))=x\lor(y\land \lnot
  z)=g_{1011}(x,y,z)$ (cf.\ Lemma~\ref{min-3aryA}).

  $\alpha$ equals $\alpha_{0010}(x,y)=x\land \lnot y$:
  \newline
  Then
  $f(x,y,z)=\alpha(x,\alpha(z,y))=x\land(y\lor \lnot
  z)=g_{0010}(x,y,z)$ (cf.\ Lemma~\ref{min-3aryA}).

  $\alpha$ equals $\alpha_{1110}(x,y)=\lnot x\lor \lnot y$:
  \newline
  Then
  $\alpha(x,\alpha(y,x))=\lnot x\lor y=\alpha_{1101}(x,y)$,
  and this $\pm$\hyp{}operation was already treated above.

  $\alpha$ equals $\alpha_{1000}(x,y)=\lnot x\land \lnot y$:
  \newline
  Then
  $\alpha(\alpha(x,y),y)= x\land \lnot y=\alpha_{0010}(x,y)$,
  and this $\pm$\hyp{}operation was already treated above.
\end{proof}


\section{Self-dual $\pm$\hyp{}preclones}\label{sec:selfdual}

In this section we consider self-dual $\pm$\hyp{}preclones. Recall that a
$\pm$\hyp{}preclone $F$ is self-dual if $F=F^{\partial}$. The self-dual $\pm$\hyp{}preclones form a
sublattice $\pmLsd_{2}$ of $\pmL_{2}$ because the meet and the join of self-dual
$\pm$\hyp{}preclones are self-dual---the former fact is trivial and the latter follows from the fact that
$\pmSg{F}^{\partial}=\pmSg{F^{\partial}}$ because $f\mapsto
f^{\partial}$ is an ``inner automorphism'' (see~\refI{Remark 5.7} for
more details).

\begin{proposition}\label{prop:selfdual1}
  Let $F\in\pmL_{2}$ be an arbitrary Boolean $\pm$\hyp{}preclone. Then
  $F\cap F^{\partial}$ and $\pmSg{F\cup F^{\partial}}$ are self-dual
  $\pm$\hyp{}preclones. We have $F\cap F^{\partial}\subseteq F\subseteq
  \pmSg{F\cup F^{\partial}}$, where each inclusion holds as an equality if and only if $F$
  is self-dual.
\end{proposition}

\begin{proof}
Since $f^{\partial\partial}=f$ we get $(F\cap
F^{\partial})^{\partial}=F^{\partial}\cap F$ and $\pmSg{F\cup
  F^{\partial}}^{\partial}=\pmSg{(F\cup
  F^{\partial})^{\partial}}=\pmSg{F^{\partial}\cup F}$.
  The claimed subset inclusions and the biconditional statement about equalities are obvious.
\end{proof}

\begin{proposition}\label{prop:selfdual2}
  \leavevmode
  \begin{enumerate}[label=\textup{(\alph*)}]

  \item\label{selfdual2a} A nontrivial self-dual $\pm$\hyp{}preclone $H$ is
    maximal in $\pmLsd_{2}$ if and only if there
    is a maximal $\pm$\hyp{}preclone $F$ such that $H=F\cap F^{\partial}$.

    There are exactly eight maximal self-dual $\pm$\hyp{}preclones, namely
    \textup{(}for the notation, see Remark~\ref{Max-remark} and Theorem~\ref{Max-theorem}\textup{)}:
\begin{align*}
&
\pmPol(\sigma_{0},\sigma_{0})\cap\pmPol(\sigma_{1},\sigma_{1}),
\quad
\pmPol(\sigma_{i},\sigma_{i}) \text{ for $i\in\{2,3,4\}$,}
\\ &
\pmPol(\leq,\geq),
\quad
\pmPol(A,\emptyset),
\quad
\pmPol(\nabla_{A},\Delta_{A}),
\quad
\pmPol(\{0\},\{1\}).
\end{align*}

  \item\label{selfdual2b} A nontrivial self-dual $\pm$\hyp{}preclone $H$ is
    minimal in $\pmLsd_{2}$ only if there
    is a minimal $\pm$\hyp{}preclone $F$ such that
    $H=\pmSg{F\cup F^{\partial}}$.

    There are exactly nine minimal self-dual $\pm$\hyp{}preclones, namely
    \textup{(}for the notation, see Theorem~\ref{min-theorem}\textup{)}:
\begin{align*}
    &
    \pmSg{f_{\mathrm{maj}}},
      \quad
    \pmSg{f_{\oplus}},
      \quad
    \pmSg{h_{1100}},
      \quad
    \pmSg{g_{0101}},
      \quad
    \pmSg{g_{1010}},
    \\ &
    \pmSg{\alpha_{0101}},
      \quad
    \pmSg{\alpha_{1010}},
      \quad
    \pmSg{\alpha_{1100}},
      \quad
    \pmSg{\alpha_{0000},\alpha_{1111}}.
\end{align*}
  \end{enumerate}

\end{proposition}

\begin{proof}

  \ref{selfdual2a} ``$\Rightarrow$'': Let $H$ be a maximal self-dual $\pm$\hyp{}preclone and let
  $F$ be a maximal $\pm$\hyp{}preclone containing $H$ ($F$ exists due to
  Proposition~\refI{5.1(3)}). Then $H':=F\cap F^{\partial}$ is not the full $\pm$\hyp{}preclone $\pmOp(A)$,
  it is self-dual (by Proposition~\ref{prop:selfdual1}) and contains $H$. By the maximality of $H$, we
  have $H=H'$.

  ``$\Leftarrow$'':
  We need to show that for each maximal $\pm$\hyp{}preclone $F$, the $\pm$\hyp{}preclone $F \cap F^{\partial}$ is maximal in $\pmLsd_{2}$.
  The eight listed $\pm$\hyp{}preclones are just all
  possible intersections $F \cap F^{\partial}$ for $F$ being one of the
  nine maximal $\pm$\hyp{}preclones (as enumerated in Theorem~\ref{Max-theorem}).
  Note that the cases with $F = \pmPol(\sigma_{0},\sigma_{0})$ and $F = \pmPol(\sigma_{1},\sigma_{1})$
  give rise to the same self-dual $\pm$\hyp{}preclone $F \cap F^{\partial}$, because $(\sigma_{0},\sigma_{0})^{\partial} = (\sigma_{1},\sigma_{1})$.

  For all the remaining seven maximal $\pm$\hyp{}preclones $F$, it holds that $F = F^{\partial}$.
  (For the cases of Theorem \ref{Max-theorem}\ref{G1a} with $i \in \{2, 3, 4\}$, this is clear because $\pmPol(\sigma_{i},\sigma_{i})$ corresponds to the self-dual clone $\Pol\sigma_{i}$.
  In all other cases, the defining relations $\rho = (\sigma,\sigma')$ satisfy $\pmPol \rho = \pmPol \rho^{\partial}$
  (since $\rho = \rho^{\partial}$ for the cases of Theorem \ref{Max-theorem}\ref{G1c},\ref{G1d},
  and $(\sigma,\sigma')^{\partial}=(\sigma',\sigma)$ for the cases of Theorem \ref{Max-theorem}\ref{G1b},\ref{G1e},
  and we can apply Proposition~\refI{5.11(3)} and Definition~\ref{rel-rules}\ref{rel-R7}).)
  Consequently, $F \cap F^{\partial} = F$, and because $F$ is a maximal $\pm$\hyp{}preclone, $F$ is obviously maximal in $\pmLsd_{2}$.

  It remains to show that also $\pmPol(\sigma_{0},\sigma_{0}) \cap \pmPol(\sigma_{1},\sigma_{1})$ is maximal in $\pmLsd_{2}$.
  To this end, it suffices to show that this is not included in any maximal $\pm$\hyp{}preclone except $\pmPol(\sigma_{0},\sigma_{0})$ and $\pmPol(\sigma_{1},\sigma_{1})$.
  Because $\pmPol(\sigma,\sigma)$ comprises all $\pm$\hyp{}operations $f$ such that $\mathring{f} \in \Pol(\sigma)$,
  it is clear from Post's lattice that
  $\pmPol(\sigma_{0},\sigma_{0}) \cap \pmPol(\sigma_{1},\sigma_{1})$
  is not included in
  $\pmPol(\sigma_i, \sigma_i)$ with $i \in \{2, 3, 4\}$.

For the remaining non\hyp{}inclusions, we observe that $\pmPol(\sigma_{0},\sigma_{0}) \cap \pmPol(\sigma_{1},\sigma_{1})$ contains every nontrivial projection (a projection whose essential argument has sign $-$), but it is straightforward to verify that none of
    $\pmPol(\leq,\geq)$,
    $\pmPol(A,\emptyset)$,
    $\pmPol(\nabla_{A},\Delta_{A})$,
    $\pmPol(\{0\},\{1\})$
contains a nontrivial projection.

  \ref{selfdual2b}:
  Let $H$ be a minimal self-dual $\pm$\hyp{}preclone and
  let $F$ be a minimal $\pm$\hyp{}preclone contained in $H$ ($F$ exists due
  to \refI{Proposition 5.1(3)}). Then $H':=\pmSg{F\cup F^{\partial}}$ is
  nontrivial, self-dual (by Proposition~\ref{prop:selfdual1}) and contained in $H$. By
  the minimality of $H$, we have $H=H'$.

We now identify those $\pm$\hyp{}preclones that are minimal in $\pmLsd_{2}$.
There are exactly fourteen $\pm$\hyp{}preclones of the form $\pmSg{F \cup F^{\partial}}$ where $F$ is one of the twenty minimal $\pm$\hyp{}preclones (as enumerated in Theorem~\ref{min-theorem}):

    $\pmSg{f_{\mathrm{maj}}}$, $\pmSg{f_{\oplus}}$,
    $\pmSg{h_{1100}}$,
    $\pmSg{g_{0101}}$, $\pmSg{g_{1010}}$,
    $\pmSg{\alpha_{0101}}$, $\pmSg{\alpha_{1010}}$,
    $\pmSg{\alpha_{1100}}$,

    $\pmSg{h_{0000},h_{1111}}$, $\pmSg{g_{0000},g_{1111}}$,
    $\pmSg{g_{0111},g_{0001}}$, $\pmSg{g_{1011},g_{0010}}$,

    $\pmSg{\alpha_{0000},\alpha_{1111}}$,
    $\pmSg{\alpha_{0001},\alpha_{0111}}$.

The first eight listed $\pm$\hyp{}preclones with one generator $g$ are the self-dual ones among the minimal $\pm$\hyp{}preclones ($g^{\partial} = g$); hence, they are obviously minimal in $\pmLsd_{2}$.
The remaining ones with two generators $g$, $g'$ satisfy $g' = g^{\partial}$, and therefore they are self-dual (by Proposition \ref{prop:selfdual1}) and of the desired form $\pmSg{F \cup F^{\partial}}$.

However, the ones with two generators except
$\pmSg{\alpha_{0000},\alpha_{1111}}$ are not minimal, because each
contains some of the self-dual $\pm$\hyp{}preclones with one generator.

In fact, we have $h_{1111}(h_{0000}(x,y,z,u),y,u,z) = h_{1100}(x,y,z,u)$
(this can be verified easily, e.g., just using the definition of $h_\alpha$).
Thus, $\pmSg{h_{1100}} \subseteq \pmSg{h_{0000},h_{1111}}$;
hence $\pmSg{h_{0000},h_{1111}}$ is not minimal in $\pmLsd_{2}$.

Analogously, $\pmSg{g_{0111}, g_{0001}}$, $\pmSg{g_{1011}, g_{0010}}$,
$\pmSg{g_{0000}, g_{1111}}$ and $\pmSg{g_{1011}, g_{0010}}$ are not
minimal in $\pmLsd_{2}$
because of the following equalities, the verification of which is
facilitated by employing the terms representing the generators of
minimal $\pm$\hyp{}preclones provided in Theorem~\ref{min-theorem}: 
\begin{align*}
& g_{0111}(g_{0001}(x,y,z),y,z)
  = ( x \wedge ( y \vee z ) ) \vee ( y \wedge z ) \\
& = ( x \wedge y ) \vee ( y \wedge z ) \vee ( z \wedge x)
  = g_{0101}(x,y,z), \\
& g_{1011}(g_{0010}(x,y,z),y,z)
  = ( x \wedge ( y \vee \neg z ) ) \vee ( y \wedge \neg z ) \\
& = ( x \wedge y ) \vee ( y \wedge \neg z ) \vee ( \neg z \wedge x)
  = g_{1010}(x,y,z), \\
& g_{1111}(g_{1111}(g_{0000}(x,y,u), g_{0000}(y,z,u), u), g_{0000}(z,x,u), u) \\
& = ( ( x \wedge y ) \vee ( y \wedge z ) ) \vee ( z \wedge x )
  = f_{\mathrm{maj}}(x,y,z,u), \\
& \alpha_{0111}(\alpha_{0001}(x,y),y)
  = ( x \wedge y ) \vee y
  = y
  = \alpha_{0101}(x,y).
\end{align*}
Finally, we see that $\pmSg{\alpha_{0000},\alpha_{1111}}$ is minimal
in $\pmLsd_{2}$, because it contains only trivial projections and
constant functions, so it cannot include any of the first eight
$\pm$\hyp{}preclones with one generator. 
\end{proof}

The knowledge of the minimal and maximal $\pm$\hyp{}preclones in the
lattice $\pmLsd_{2}$ according to Proposition~\ref{prop:selfdual2} is of special interest because $\pmLsd_{2}$ is
atomic and coatomic. In particular, we get a completeness criterion
which slightly specializes Theorem~\ref{thm:complete1}.

\begin{proposition}\label{prop:SDlattice}
  The lattice $\pmLsd_{2}$ is atomic and coatomic, i.e., each
    self-dual $\pm$\hyp{}preclone is contained in a maximal self-dual
    $\pm$\hyp{}preclone and contains a minimal self-dual $\pm$\hyp{}preclone.
\end{proposition}

\begin{proof}
  Let $H=H^{\partial}$ be a self-dual $\pm$\hyp{}preclone. Then there exists
  a maximal $\pm$\hyp{}preclone $F$ containing $H$ (since $\pmL_{2}$ is
  coatomic), thus $F^{\partial}$ also
  contains $H^{\partial}=H$, i.e., $H\subseteq F\cap F^{\partial}$,
  the latter being a maximal self-dual $\pm$\hyp{}preclone because of
  Proposition~\ref{prop:selfdual2}\ref{selfdual2a}. Analogously, $H$ contains a
  minimal self-dual $\pm$\hyp{}preclone $\pmSg{F\cup F^{\partial}}$ for some
  minimal $\pm$\hyp{}preclone $F$.
\end{proof}

\begin{theorem}[Completeness theorem for self-dual $\pm$\hyp{}preclones]
\label{thm:complete2}
A set $F\subseteq\pmOp(A)$ with the property $F=F^{\partial}$
     is complete \textup{(}i.e., it
     generates the full \textup{(}self-dual\textup{)} $\pm$\hyp{}preclone $\pmOp(A)$\textup{)}
    if and only if for each of the following eight $\pm$\hyp{}relations there
    exists a $\pm$\hyp{}operation in $F$ which does not $\pm$\hyp{}preserve this
    relation:
    $(\sigma_{i},\sigma_{i})$ for $i\in\{0,2,3,4\}$,
    $(\leq,\geq)$, $(A,\emptyset)$,
    $(\nabla_{A},\Delta_{A})$, $(\{0\},\{1\})$.

\end{theorem}

\begin{proof} We can apply Theorem~\ref{thm:complete1}. Thus it
  remains to show that $(\sigma_{1},\sigma_{1})$ need not appear in
  the listed $\pm$\hyp{}relations. This follows from the fact that for a
  $\pm$\hyp{}operation $g\in F$ and its dual $f:=g^{\partial}\in
  F^{\partial}=F$ the following holds: $g$ does not $\pm$\hyp{}preserve
  $(\sigma_{0},\sigma_{0})$ if and only if $f$ does not $\pm$\hyp{}preserve
  $(\sigma_{1},\sigma_{1})$ (cf.\ Proposition~\refI{5.11(i)} with
  $h:=\partial$).
\end{proof}


\section{The lattice $\pmL_{2}$ of Boolean $\pm$\hyp{}preclones}
\label{sec:latticeL2}

Boolean $\pm$\hyp{}preclones generalize Boolean clones and they are special
 preclones (for this notion see Remark~\refI{2.4(D)}). It is well known that there are countably many Boolean clones
 (Post's lattice). Although $\pm$\hyp{}preclones are relatively close
 (among all preclones)  to clones we shall see now that the lattice
 $\pmL_{2}$ is already uncountable. 
This result is due to Vojt\v{e}ch David (personal communication, Sept.\ 2026)

\begin{theorem}[V. David]\label{thm:uncountable}
  There are uncountably many Boolean $\pm$\hyp{}preclones.
\end{theorem}

Before we can start the proof we need some notions, notations and
facts concerning clonoids.

For sets $K, L \subseteq \Op(A)$, the \New{composition} of $K$ and $L$ is the set
\[
\begin{split}
KL :=
\{ \, h \in Op(A) \mid {} & \exists m \in \N_{+} \, \exists n \in \N_{+} \, \exists f \in K^{(n)} \, \exists g_1, \dots, g_n \in L^{(m)} \colon
\\ &
h = f(g_1, \dots, g_n) \, \}
\end{split}
\]
of all possible compositions of operations where the outer function comes from $K$ and all inner functions come from $L$.
Using this notion of function class composition, we can easily discuss further algebraic concepts.
For example, a clone on $A$ is a set $C \subseteq \Op(A)$ that contains all projections and $C C \subseteq C$.

Let now $C_1$ and $C_2$ be two clones on $A$, referred to as the \New{source clone} and the \New{target clone}, respectively.
A set $K \subseteq \Op(A)$ is called a \New{$(C_1,C_2)$-clonoid} if $K C_1 \subseteq K$ and $C_2 K \subseteq K$.
The set of all $(C_1,C_2)$-clonoids constitutes a closure system, i.e., $\Op(A)$ is a $(C_1,C_2)$-clonoid and arbitrary intersections of $(C_1,C_2)$-clonoids are again $(C_1,C_2)$-clonoids.

We need the following two clones on $A = \{0,1\}$.
\begin{itemize}
\item
  $\clUinf := \{ \, f \in \Op(\{0,1\}) \mid \bigwedge f^{-1}(1) \neq
  (0, \dots, 0) \, \}$ (here, the meet is taken with respect to the
  componentwise ordering of tuples, with $0 < 1$),
\item $\clLambda$ is the clone of conjunctions and
  constants.
\end{itemize}

\begin{theorem}[{\cite[Theorem 5.5]{Leh2026b}, see also \cite[Theorem
    5.3, Table 2]{Leh2026a}}]
\label{thm:UinfLambda}
There are an uncountable infinitude of $(\clUinf, \clLambda)$-clonoids.
\end{theorem}

We are going to prove that the poset of $(\clUinf, \clLambda)$-clonoids embeds into the poset of Boolean $\pm$\hyp{}preclones.
To this end, we need the following two maps.

\begin{definition}
\label{def:sharp-flat}
For $f \in \Op(A)$ of arity $n$, let $f^\sharp$ be the $(n+1)$\hyp{}ary $\pm$\hyp{}operation
\[
f^\sharp(x^{+}, y_1^{-}, \dots, y_n^{-}) := x^{+} \wedge f(y_1^{-}, \dots, y_n^{-}),
\]
and for $C \subseteq \Op(A)$, let $C^\sharp$ be the $\pm$\hyp{}preclone generated by $\{ \, f^\sharp \mid f \in C \,\}$.
Conversely, for $P \subseteq \pmOp(A)$, let $P^\flat := \{ \, f \in \Op(A) \mid f^\sharp \in P \, \}$.
\end{definition}

\begin{proposition}[V.\ David]
\label{prop:uncountable:aux}
If $C$ is a nonempty $(\clUinf, \clLambda)$-clonoid, then $C^\sharp$ is a Boolean $\pm$\hyp{}preclone and $(C^\sharp)^\flat = C$.
\end{proposition}

\begin{proof}
The inclusion $C \subseteq (C^\sharp)^\flat$ is immediate from Definition~\ref{def:sharp-flat}.

For the converse inclusion, let $D$ be the set of $\pm$\hyp{}operations $h$, with variables $u_1, \dots, u_n$, that can be expressed as
\begin{equation}
h(u_1, \dots, u_n) = x \wedge g ( z_1 \wedge \varphi_1(u_1, \dots, u_n), \dots, z_m \wedge \varphi_m(u_1, \dots, u_n) ),
\label{eq:D}
\end{equation}
where
\begin{itemize}
\item $g$ is an $m$\hyp{}ary operation from $C$, $m \geq 1$,
\item $x \in \{u_1, \dots, u_n\}$ is a positive variable,
\item each $z_j$ is either constant $0$ or a negative variable from $\{u_1, \dots, u_n\}$,
\item $\varphi_1, \dots, \varphi_m$ are arbitrary $n$\hyp{}ary $\pm$\hyp{}operations of suitable signa.
\end{itemize}

\textit{Claim 1.}
If $f \in C$, then $f^\sharp \in D$.

\textit{Proof of Claim 1.}
We can write $f^\sharp$ as
\[
f^\sharp(x^{+}, y_1^{-}, \dots, y_n^{-})
= x^{+} \wedge f(y_1^{-}, \dots, y_n^{-})
= x^{+} \wedge f(y_1^{-} \wedge y_1^{-}, \dots, y_n^{-} \wedge y_n^{-}),
\]
which is clearly of the form \eqref{eq:D}.
$\diamond$

\textit{Claim 2.}
$\id_A$ (with signum $+$) is in $D$.

\textit{Proof of Claim 2.}
We can write $\id_A$ as
\[
\id_A(x^{+}) = x^{+} \wedge 1(0 \wedge 0),
\]
which is of the form \eqref{eq:D}.
Note that the unary constant function $1$ belongs to the
$(\clUinf,\clLambda)$-clonoid $C$ because $1 \in \clLambda$ and hence
$1 \in \clLambda C \subseteq C$ (note $C\neq\emptyset$).
$\diamond$

\textit{Claim 3.}
$D$ is closed under composition.

\textit{Proof of Claim 3.}
Let $f \in D^{(n)}$ and $g_1, \dots, g_n \in D^{(m)}$ with $\sgn(f) = (s_1, \dots, s_n) \in \{\mathord{+}, \mathord{-}\}^n$ and for some $(s'_1, \dots, s'_m) \in \{\mathord{+}, \mathord{-}\}^m$, $\sgn(g_i) = (s'_1, \dots, s'_m)$ if $s_i = \mathord{+}$ and $\sgn(g_i) =  (\mathord{-} \cdot s'_1, \dots, \mathord{-} \cdot s'_m)$ if $s_i = \mathord{-}$, for $i \in \{1, \dots, m\}$.
Without loss of generality, $s_n = \mathord{+}$ and $s_1 = \dots = s_\ell = \mathord{-}$ for some $\ell < n$ and
\[
f(x_1, \dots, x_n) = x_n \wedge \gamma ( x_1 \wedge \varphi_1(x_1, \dots, x_n), \dots, x_\ell \wedge \varphi_\ell(x_1, \dots, x_n), 0, \dots, 0)
\]
for some $\gamma \in C$ of arity $r$, with $\ell \leq r < n$ and for some $\pm$\hyp{}operations $\varphi_1, \dots, \varphi_\ell$ of suitable signa,
and for $i \in \{1, \dots, n\}$,
\begin{multline*}
g_i(x_1, \dots, x_m) = \\
x_{j_i} \wedge \gamma_i ( x_{j_{i,1}} \wedge \rho_{i,1}(x_1, \dots, x_m), \dots, x_{j_{i,\ell_i}} \wedge \rho_{i,\ell_i}(x_1, \dots, x_m), 0, \dots, 0)
\end{multline*}
for some $\gamma_i \in C$ of arity $r_i$, with $\ell_i \leq r_i < m$, for some $j_i, j_{i,1}, \dots, j_{i,\ell_i} \in \{1, \dots, m\}$ such that the signum of $x_{j_i}$ is $+$ in $g_i$ and the signum of each of $x_{j_{i,1}}, \dots, x_{j_{{i,\ell_i}}}$ is $-$, and for some $\pm$\hyp{}operations $\rho_{i,1}, \dots, \rho_{i,\ell_i}$ of suitable signa,
Then
\[
\begin{split}
& f(g_1, \dots, g_n)(x_1, \dots, x_m)
\\ &
= g_n(x_1, \dots, x_m)
\\ &
\phantom{{}={}}
\wedge
\gamma (
g_1(x_1, \dots, x_m) \wedge \varphi_1(g_1, \dots, g_n)(x_1, \dots, x_m),
\dots,
\\ & \phantom{{} = {} \wedge \gamma (}
g_\ell(x_1, \dots, x_m) \wedge \varphi_\ell(g_1, \dots, g_n)(x_1, \dots, x_m),
0, \dots, 0
)
\\ &
=
x_{j_n} \wedge \gamma_n ( x_{j_{n,1}} \wedge \rho_{n,1}(x_1, \dots, x_m), \dots, x_{j_{n,\ell_n}} \wedge \rho_{n,\ell_n}(x_1, \dots, x_m), 0, \dots, 0)
\\ &
\phantom{{}={}}
\wedge
\gamma (
g_1(x_1, \dots, x_m) \wedge \varphi_1(g_1, \dots, g_n)(x_1, \dots, x_m),
\dots,
\\ & \phantom{{}={} \wedge \gamma (}
g_\ell(x_1, \dots, x_m) \wedge \varphi_\ell(g_1, \dots, g_n)(x_1, \dots, x_m),
0, \dots, 0
)
\\ &
=
x_{j_n} \wedge
\hat{\gamma}(
x_{j_n} \wedge \gamma_n ( x_{j_{n,1}} \wedge \rho_{n,1}(x_1, \dots, x_m),
\dots,
x_{j_{n,\ell_n}} \wedge \rho_{n,\ell_n}(x_1, \dots, x_m),
\\ & \phantom{{}={} x_{j_n} \wedge \gamma(}
0, \dots, 0,
x_{j_1} \wedge \zeta_1 (x_1, \dots, x_m),
\dots,
x_{j_\ell} \wedge \zeta_\ell (x_1, \dots, x_m),
0, \dots, 0
)
,
\end{split}
\]
where
\[
\hat{\gamma}(x_1, \dots, x_{r_n}, y_1, \dots, y_r)
= \gamma_n(x_1, \dots, x_{r_n}) \wedge \gamma(y_1, \dots, y_r),
\]
and for $i \in \{1, \dots, \ell\}$,
\[
\begin{split}
\zeta_i (x_1, \dots, x_m) = {}
& \gamma_i ( x_{j_{i,1}} \wedge \rho_{i,1}(x_1, \dots, x_m), \dots, x_{j_{i,\ell_i}} \wedge \rho_{i,\ell_i}(x_1, \dots, x_m), 0, \dots, 0)
\\ & {}
\wedge
\varphi_i (g_1, \dots, g_n) (x_1, \dots, x_m)
\end{split}
\]
Because $C$ is closed under minors and $\Lambda C \subseteq C$, it follows that $\hat{\gamma} \in C$.
Therefore, $f(g_1, \dots, g_n) \in D$.
$\diamond$

\textit{Claim 4.}
$D^\flat \subseteq C$.

\textit{Proof of Claim 4.}
Let $q \in D^\flat$ be $n$\hyp{}ary, that is,
the $(n+1)$\hyp{}ary $\pm$\hyp{}operation
\[
q^\sharp(x^{+}, y_1^{-}, \dots, y_n^{-}) = x \wedge q(y_1, \dots, y_n)
\]
is an element of $D$.
Note that since $x$ is the only positive variable in $q^\sharp$, in any expression of the form \eqref{eq:D} representing $q^\sharp$, the first variable (on the left side of $\wedge$) must be $x$.
Fix one such expression; thus, we have, for some $g \in C$ and for some $\pm$\hyp{}operations $\varphi_1, \dots, \varphi_m$,
\[
x \wedge q(y_1, \dots, y_n)
=
x \wedge g ( z_1 \wedge \varphi_1 (x, y_1, \dots, y_n), \dots, z_m \wedge \varphi_m(x, y_1, \dots, y_n) ),
\]
where each $z_i$ is either constant $0$ or one of the negative variables $y_1, \dots, y_n$.

Substituting $1$ for $x$ on both sides, we obtain
\[
q(y_1, \dots, y_n)
=
g ( z_1 \wedge \varphi_1 (1, y_1, \dots, y_n), \dots, z_m \wedge \varphi_m(1, y_1, \dots, y_n) ),
\]
Each of the $n$\hyp{}ary operations $z_j \wedge \varphi_j (1, y_1, \dots, y_n)$ is either constant $0$ or has the variable $z_j$ equal to $1$ in each of its satisfying assignments;
in both cases it belongs to $\clUinf$.
Hence $q \in C \clUinf \subseteq C$.
$\diamond$

By Claims 2 and 3, $D$ is a $\pm$\hyp{}preclone, and together with Claim 1, we get $C^\sharp \subseteq D$.
Claim 4 then gives $(C^\sharp)^\flat \subseteq D^\flat \subseteq C$, as desired.
\end{proof}

\begin{proof}[Proof of Theorem~\ref{thm:uncountable}]
Proposition~\ref{prop:uncountable:aux} shows that the mapping $C \mapsto C^\sharp$ is an embedding of the poset of nonempty $(\clUinf, \clLambda)$-clonoids into the poset of Boolean $\pm$\hyp{}preclones.
Since the lattice of $(\clUinf, \clLambda)$-clonoids is uncountable by Theorem~\ref{thm:UinfLambda}, the result follows.
\end{proof}

Now we develop some tools to investigate parts of the lattice
$\pmL_{2}$.

\begin{definition}\label{def:preimage}
  Let $F\leq \Op(A)$ be a Boolean clone, i.e., $F\in\cL_{2}$. A
  $\pm$\hyp{}preclone $P\leq\pmOp(A)$ is called a \New{preimage of $F$} if
  $\Sg{\mathring P}=F$, i.e., taking all operations in $P$ and
  ignoring signa we can generate all operations in $F$.  Let $\pF$
  denote the set of all preimages of $F$, called 
  \New{preimage class of $F$}.
\end{definition}

A description of $\Sg{\mathring P}$
by invariant relations was already given in Proposition~\ref{F1}.
One can see that $\Sg{\mathring P}$ can be obtained from
  $\mathring P$ by 
  adding fictitious variables and identification of variables only.

Clearly, $\{\pF\mid F\in \cL_{2}\}$ is a partition of
$\pmL_{2}$ into preimage classes. Therefore we ask for the structure of these
preimage classes.

\begin{proposition}\label{prop:preimage1} Let $F\in\cL_{2}$ and let
  $G$ be a generating system of the Boolean clone $F$.
  Each preimage class $\pF$ is a semi-interval in $\pmL_{2}$: it has a greatest element, namely the
  $\pm$\hyp{}preclone
  \begin{align*}
    P_{F}:=\{f\in\pmOp(A)\mid \mathring f\in
    F\}=\pmSg{\{\id_{A}^{-}\}\cup\{g^{(+,\ldots,+)}\mid g\in G\}}
  \end{align*}
but, in general, there may exist several minimal elements in $\pF$.
\end{proposition}

\begin{figure}
  \begin{center}
  \includegraphics{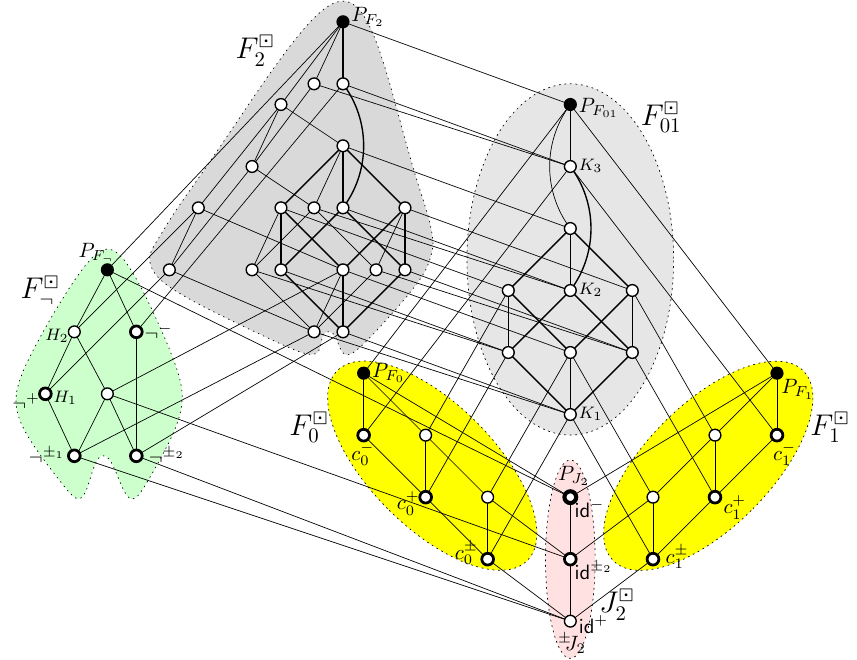}
  \end{center}
\caption{The lattice of the unary $\pm$\hyp{}preclones with its six 
preimage classes $\pF[J_2]$, $\pF_{0}$, $\pF_{1}$, $\pF[F_{01}]$,
$\pF[F_{\neg}]$ and $\pF[F_{2}]$ (the greatest elements in each
preimage class as well as the join\hyp{}irreducibles are
marked)\label{fig:preimagesA}} 

\end{figure}

\begin{proof} It is straightforward to see that the left side $P_{F}$ is a
  $\pm$\hyp{}preclone which obviously is the greatest $\pm$\hyp{}preclone in
  $\pF$. The indicated generating system generates $P_{F}$ since one
  can change the signum of each variable arbitrarily by composition
  with $\id_{A}^{-}$. 

 Moreover, the interval property is clear since
  for $\pm$\hyp{}preclones $P\leq P'$ with $\Sg{\mathring P}=\Sg{\mathring
    P'}=F$ we obviously have $\Sg{\mathring P''}=F$ for each
  $\pm$\hyp{}preclone $P\leq  P''\leq P'$. 

$\pF$ is not always an interval as the example $F=\Sg{\neg}$ shows
(see Figure~\ref{fig:preimagesA}).
\end{proof}

According to Theorem~\refI{4.1} each $\pm$\hyp{}preclone can be characterized by
invariant relations $\rho=(\rho_{+},\rho_{-})$. The next proposition
gives a relational characterization of the preimage classes $\pF$.

\begin{proposition}\label{prop:relpreimages}
 Let $F\in\cL_{2}$ be a Boolean clone. Then a $\pm$\hyp{}preclone
 $H\in\pmL_{2}$ belongs to $\pF$ if and only if 
  $\Inv F=\{\sigma\in\Rel(A)\mid H\pmpreserves(\sigma,\sigma)\}$.
\end{proposition}
\begin{proof} 
Note $H\in\pF\iff\Sg{\mathring H}=F\iff\Inv\Sg{\mathring H}=\Inv F$.
By Proposition~\refI{3.10} we have $\Sg{\mathring
    H}=\Pol\{\sigma\in\Rel(A)\mid (\sigma,\sigma)\in\pmInv H\}$.
Consequently, $\Inv\Sg{\mathring H}=\Inv\Pol\{\sigma\in\Rel(A)\mid
(\sigma,\sigma)\in\pmInv H\}=\{\sigma\in\Rel(A)\mid
(\sigma,\sigma)\in\pmInv H\}$, which completes the proof (the latter
equality follows from the 
fact that the set $\{\sigma\in\Rel(A)\mid
(\sigma,\sigma)\in\pmInv H\}$ is a relational clone (cf.\ Definition~\refI{2.7}),
which immediately follows from Definition~\refI{2.11} for the
$S$\hyp{}relational clone $H$). 
\end{proof}

We are going to describe $\pF$
for special Boolean clones $F$.

At first we deal with \New{unary Boolean $\pm$\hyp{}preclones} $P$, i.e.,
$\pm$\hyp{}preclones where
$\mathring{f}$ is essentially unary for each $f\in P$. Such
$\pm$\hyp{}preclones are completely characterized by their 
unary or binary $\pm$\hyp{}operations as Lemma~\ref{lem:unarypreclones}
will show.

For a unary operation $g\in\Opa[1](A)$ we consider the following
associated $\pm$\hyp{}operations ($x,y\in A$):
\begin{align*}
  g^{+}(x)&:=g(x),\; \sgn(g^{+})=(+),&g^{\pm_{^{1}}}(x,y)&:=g(x),\;\sgn(g^{\pm_{1}})=(+,-), \\
 g^{-}(x)&:=g(x),\; \sgn(g^{-})=(-), & g^{\pm_{2}}(x,y)&:=g(y),\; \sgn(g^{\pm_{2}})=(+,-).
\end{align*}
For $c_{i}$ being the unary
Boolean operation constant $i$ (for $i\in\{0,1\}$) we have $c_{i}^{\pm_{1}}=c_{i}^{\pm_{2}}$. Therefore we
denote this $\pm$\hyp{}operation by $c_{i}^{\pm}$ ($i\in\{0,1\}$). 
Moreover, for short we write $\id$ instead of $\id_{A}$.

Note that all generators of unary minimal $\pm$\hyp{}preclones are among
the above $\pm$\hyp{}operations (for the notation, cf.\ Theorem~\ref{min-theorem}\ref{min-thmD}):
$\alpha_{0000}=c_{0}^{\pm}$, $\alpha_{1111}=c_{1}^{\pm}$,
$\alpha_{0101}=\id^{\pm_{2}}$, $\alpha_{1100}=\neg^{\pm_{1}}$,
$\alpha_{0101}=\neg^{\pm_{2}}$. These generate the minimal elements in the
lattice of unary Boolean $\pm$\hyp{}preclones in Figure~\ref{fig:preimagesA}.

\begin{lemma}\label{lem:unarypreclones} Let $P$ be a unary Boolean
  $\pm$\hyp{}preclone. For $f\in P$ let $f'$ be the $\pm$\hyp{}operation
  obtained from $f$ by identifying all arguments of the same
  sign. Then $\pmSg{f'}=\pmSg{f}$.
\end{lemma}
Note that $f'$ must have one of the following signa: $(+)$, $(-)$,
$(+,-)$, $(-,+)$ but by permutation of arguments we may exclude
$(-,+)$.

\begin{proof}
  $P$ is a unary $\pm$\hyp{}preclone, so $f$ has at most one essential
  argument. This immediately implies the statement of the lemma.
\end{proof}

Therefore, in order to find all unary $\pm$\hyp{}preclones, it is enough to
find all $\pm$\hyp{}preclones generated by $\pm$\hyp{}operations of the form
$g^{+},g^{-},g^{\pm_{1}},g^{\pm_{2}}$ for some $g\in\Opa[1](A)$. The result is given in Theorem~\ref{thm:preimage3}.

\begin{lemma}\label{lem:preimage2} Let $g\in\Opa[1](A)$. Then
  \begin{enumerate}[label={\upshape(\alph*)}]
  \item \label{preimage2-a}
   $\pmSg{g^{\pm_{1}}}\subsetneqq \pmSg{g^{+}}$ \textup{(}for $g\neq\id$\textup{)},
  \item \label{preimage2-b} $\pmSg{g^{\pm_{2}}}\subsetneqq
    \pmSg{g^{-}}$,

\item \label{preimage2-c}
  $\pmSg{g^{-}}\subsetneqq\pmSg{g^{\pm_{1}},\id^{-}}$ \textup{(}for $g\neq\id$\textup{)}.

\end{enumerate}

\end{lemma}

\begin{proof}
\ref{preimage2-a}: We have $\pmSg{g^{\pm_{1}}}\subseteq \pmSg{g^{+}}$
since $g^{\pm_{1}}(x,y)=x=(\nabla^{-}g^{+})(y,x)$. The inclusion is strict
since there exists $\rho\in\pmRel(A)$ such that
$g^{\pm_{1}}$ does $\pm$\hyp{}preserve $\rho$ but
$g^{+}$ does not (in fact, for $g\neq \id$ there exist $a\neq b$ such that
$g(a)=b$ and the unary $\pm$\hyp{}relation $\rho:=(\emptyset,\{a\})$ will do
the job).

\ref{preimage2-b}: We have $\pmSg{g^{\pm_{2}}}\subseteq \pmSg{g^{-}}$
since $g^{\pm_{2}}=\nabla^{+}g^{-}$. The inclusion is strict since 
$\rho\in\pmInv g^{\pm_{2}}\setminus\pmInv g^{-}$ for
$\rho=(A,\emptyset)$ (cf.\ Lemma~\ref{Max-L4}).

\ref{preimage2-c}: We have
$\pmSg{g^{-}}\subseteq\pmSg{g^{\pm_{1}},\id^{-}}$ since
$g^{-}(x)=g^{\pm_{1}}(\id^{-}(x),x)$. The inclusion is strict since
$\rho\in\pmInv g^{-}\setminus\pmInv \{g^{\pm_{1}},\id^{-}\}$ for the
binary $\pm$\hyp{}relation $\rho$ with
$\rho_{+}=\{(0,1),(g^{2}(0),g^{2}(1))\}$ and $\rho_{-}:=\{(g(0),g(1))\}$
(note that $g^{3}=g$ and $\rho_{+}\neq\rho_{-}$ if $g\neq\id$;
therefore $\id^{-}$ does not $\pm$\hyp{}preserve $\rho$, but $g^{-}$ does).
\end{proof}

There are six monoids of unary Boolean functions:
$M_{\id}:=\{\id_{A}\}$, $M_{0}:=\{\id_{A},c_{0}\}$, $M_{1}:=\{\id_{A},c_{1}\}$,
  $M_{01}:=\{\id_{A},c_{0},c_{1}\}$, $M_{\lnot}:=\{\id_{A},\lnot\}$,
  $M_{2}:=\{\id_{A},c_{0},c_{1},\lnot\}=\Opa[1](A)$. The Boolean
  clones generated by these $M_{i}$ coincide with $M_{i}$ up to
  fictitious variables and shall be denoted by $F_{i}:=\Sg{M_{i}}$ for
  $i\in\{\id,0,1,01,\lnot,2\}$, note $F_{\id} = \Sg{M_{\id}}=J_{2}$. 
 The
corresponding preimage classes are 
  described in the next theorem (each $\pm$\hyp{}preclone in each class is given by a
  generating system). 

\begin{theorem}\label{thm:preimage3} 

  The lattice of unary Boolean $\pm$\hyp{}preclones has altogether 51
  elements. There are 12 join\hyp{}irreducible elements generated by
$\id_{A}^{\pm_{2}}$, $\id_{A}^{-}$, $c_0^{\pm}$,
$c_0^{+},c_0^{-}$, $c_1^{\pm}$, $c_1^{+},c_1^{-}$, $\neg^{\pm_{1}}$,
$\neg^{\pm_{2}}$,  $\neg^{+}$,  $\neg^{-}$,
where we have the following inclusions:
  \begin{align*}
  \pmSg{\id_{A}^{\pm_{2}}}&<\pmSg{\id_{A}^{-}},&
 \pmSg{c_0^{\pm}}&< \pmSg{c_0^{+}}< \pmSg{c_0^{-}},\\
  \pmSg{\neg^{\pm_{1}}}&< \pmSg{\neg^{+}},&
  \pmSg{c_1^{\pm}}&< \pmSg{c_1^{+}}< \pmSg{c_1^{-}},\\
  \pmSg{\neg^{\pm_{2}}}&< \pmSg{\neg^{-}}.
  \end{align*}

The unary and binary operations in each $\pm$\hyp{}preclone can be read
from the lattice in Figure~\ref{fig:preimagesA}. The six preimage
classes are as follows. Each $\pm$\hyp{}preclone has a generating system of
at most four $\pm$\hyp{}operations \textup{(}it is the join of at most four join\hyp{}irreducibles\textup{)}.
  \begin{enumerate}[label={\upshape(\roman*)}]
  \item\label{preimage3-i} $\pF[J_{2}]=\{\pmSg{\id_{A}^{+}},
    \pmSg{\id_{A}^{\pm_{2}}}, \pmSg{\id_{A}^{-}}\}$, 
    note $\pmJ_{2}=\pmSg{\id_{A}^{+}}$, $P_{J_{2}}=\pmSg{\id_{A}^{-}}$, 
   \hfill $|\pF[J_{2}]|=3$,

\item\label{preimage3-ii}
  $\pF[F_{0}]=\{\pmSg{c_{0}^{\pm}},\pmSg{c_{0}^{+}},\pmSg{c_{0}^{-}},
    \pmSg{\id^{\pm_{2}},c_{0}^{\pm}},\pmSg{\id^{\pm_{2}},c_{0}^{+}},
    \pmSg{\id^{\pm_{2}},c_{0}^{-}}\}$, \hfill $|\pF[F_{0}]|=6$,

\item\label{preimage3-iii}
  $\pF[F_{1}]=\{\pmSg{c_{1}^{\pm}},\pmSg{c_{1}^{+}},\pmSg{c_{1}^{-}},
    \pmSg{\id^{\pm_{2}},c_{1}^{\pm}},\pmSg{\id^{\pm_{2}},c_{1}^{+}},\pmSg{\id^{\pm_{2}},c_{1}^{-}}\}$,\hfill $|\pF[F_{1}]|=6$,

\item \label{preimage3-iv}
 $\pF[F_{01}]=\{
      \pmSg{c_{0}^{\pm},c_{1}^{\pm}},
      \pmSg{\id^{\pm_{2}},c_{0}^{\pm_{1}},c_{1}^{\pm}},\\
    \hspace*{7ex}\pmSg{c_{0}^{\pm},c_{1}^{+}},
    \pmSg{\id^{\pm_{2}},c_{0}^{\pm},c_{1}^{+}},
    \pmSg{c_{1}^{\pm},c_{0}^{+}},
    \pmSg{\id^{\pm_{2}},c_{1}^{\pm},c_{0}^{+}},\\
    \hspace*{7ex}\pmSg{c_{0}^{+},c_{1}^{+}},
    \pmSg{\id^{\pm_{2}}, c_{0}^{+},c_{1}^{+}},
    \pmSg{c_{0}^{-},c_{1}^{-}},
    \pmSg{\id^{\pm_{2}}, c_{0}^{-},c_{1}^{-}}\}$,
\hfill $|\pF[F_{01}]|=10$,

note $\pmSg{\id^{\pm_{2}},
  c_{0}^{-},c_{1}^{-}}=P_{F_{01}}=\pmSg{\id^{-},c_{0}^{+},c_{1}^{+}}
    =P_{F_{0}}\lor P_{F_{1}}$,

\item\label{preimage3-v}
 $\pF[F_{\lnot}]=\{ \pmSg{\neg^{\pm_{1}}},\pmSg{\neg^{\pm_{2}}},
             \pmSg{\neg^{+}},\pmSg{\neg^{-}},\pmSg{\id^{\pm_{2}},\neg^{\pm_{1}}},
        \pmSg{\id^{\pm_{2}},\neg^{+}}, \pmSg{\id^{\pm_{2}},\neg^{-}}\}$,\\ 
\hspace*{\fill}$|\pF[F_{\lnot}]|=7$,

note $\pmSg{\id^{\pm_{2}},\neg^{-}}=P_{F_{\neg}}=\pmSg{\id^{-},\neg^{+}}$,

\item\label{preimage3-vi} 
$\pF[F_{2}]=\{\pmSg{\neg^{\pm_{1}},c_{0}^{\pm}}, 
              \pmSg{\neg^{\pm_{1}},c_{0}^{+}},\\
    \hspace*{7ex}\pmSg{\neg^{\pm_{1}},c_{1}^{+} }, 
               \pmSg{\neg^{\pm_{1}},c_{0}^{+},c_{1}^{+}},\\[0.7ex]
    \hspace*{7ex} \pmSg{\neg^{+},c_{0}^{\pm}},
               \hspace*{1ex}\pmSg{\neg^{+},\id^{\pm_{2}},c_{0}^{\pm}},\\
    \hspace*{7ex} \pmSg{\neg^{+},c_{0}^{+}},
              \hspace*{1ex} \pmSg{\neg^{+},\id^{\pm_{2}},c_{0}^{+}},
               \pmSg{\neg^{+},c_{0}^{-}},\\[0.7ex]
    \hspace*{7ex}\pmSg{\neg^{\pm_{2}}, c_{0}^{\pm}},
          \pmSg{\neg^{\pm_{2}}, \id^{\pm_{2}},c_{0}^{\pm_{1}}},\\
    \hspace*{7ex}\pmSg{\neg^{\pm_{2}},c_{0}^{+}},
          \pmSg{\neg^{\pm_{2}}, \id^{\pm_{2}},c_{0}^{+}},\\
    \hspace*{7ex} \pmSg{\neg^{\pm_{2}},c_{1}^{+}},
             \pmSg{\neg^{\pm_{2}}, \id^{\pm_{2}},c_{1}^{+}},\\
    \hspace*{7ex}\pmSg{\neg^{\pm_{2}}, c_{0}^{+},c_{1}^{+}},
            \pmSg{\neg^{\pm_{2}}, \id^{\pm_{2}}, c_{0}^{+},c_{1}^{+}},\\
    \hspace*{7ex}\pmSg{\neg^{\pm_{2}}, c_{0}^{-}},
            \pmSg{\neg^{\pm_{2}}, \id^{\pm_{2}}, c_{0}^{-}}
\}$ 
\hfill $|\pF[F_{2}]|=19$, 

note
$\pmSg{\neg^{\pm_{2}}, \id^{\pm_{2}}, c_{0}^{-}}=P_{F_{2}} =\pmSg{\id^{-},c_{0}^{+},c_{1}^{+},\neg^{+}}
                 =P_{F_{\neg}}\lor P_{F_{01}}$.

  \end{enumerate}
\end{theorem}

\begin{proof}[Sketch of the proof]
All necessary computations are more or less straightforward. In order
to show that $\pmSg{G_{1}}\subseteq\pmSg{G_{2}}$ (for $\pm$\hyp{}preclones
from the lists in the theorem generated by $G_{1}$ and $G_{2}$) for
each $g_{1}\in G_{1}$ one has to present a term (with the operations
from $G_{2}$) whose term operation coincides with $g_{1}$. In order to
show that $\pmSg{G_{1}}\neq\pmSg{G_{2}}$, it is enough to find a $\pm$\hyp{}relation
$\rho=(\rho_{+},\rho_{-})$ which is $\pm$\hyp{}preserved by one of $G_{1}$ or
$G_{2}$ and not by the other (e.g.,
$\rho\in\pmInv G_{1}\setminus\pmInv G_{2}$ shows that the above
inclusion is strict). 
In this way one also can show that
$\pmSg{G_{1}}\lor\pmSg{G_{2}}=\pmSg{G_{3}}$ by proving
$\pmSg{G_{1}\cup G_{2}}\subseteq\pmSg{G_{3}}$ and $\pmSg{G_{3}}\subseteq\pmSg{G_{1}\cup G_{2}}$.

Finally, for the correctness of Figure~\ref{fig:preimagesA} one has to check
that for each $\pm$\hyp{}preclone $P$ with lower covers
$P_{1},\dots,P_{n}$ we have $P=P_{i}\lor P_{j}$ for all
$i,j\in\{1,\dots,n\}$, $i\neq j$. This completes the proof since all
possible join\hyp{}irreducibles and their joins appear in
Figure~\ref{fig:preimagesA}.

We shall prove the general structure of
$\pmL_{2}$ and some selected examples but we do not give each detail
(left as an exercise that can be done, e.g., with the help of a computer program).

We start with the \emph{join\hyp{}irreducibles}:
Each join\hyp{}irreducible unary $\pm$\hyp{}preclone must be generated by a single
essentially unary $\pm$\hyp{}operation. Because there are only four unary
Boolean operations ($\id_{A},c_{0},c_{1},\neg$), the $12$
$\pm$\hyp{}operations in the Theorem list \textsl{all} possible candidates (clearly
$\id_{A}^{+}$ and $\id_{A}^{\pm_{1}}$ are excluded because they
generate the trivial $\pm$\hyp{}preclone $\pmJ_{2}$) and
because of the indicated inclusions they are all different and really
are join\hyp{}irreducible. To check an example: we have
$\pmSg{c_{0}^{+}}\subsetneqq\pmSg{c_{0}^{-}}$ since
$c_{0}^{+}(x)=c_{0}^{-}(c_{0}^{-}(x))$ and $\rho\in\pmInv
c_{0}^{+}\setminus\pmInv c_{0}^{-}$ for $\rho=(\{0\},\emptyset)$.

Now it ``only'' remains to consider all possible joins of the 12
join\hyp{}irreducible $\pm$\hyp{}preclones. Each $\pm$\hyp{}preclone belongs to one
of the six preimage classes which we discuss in more detail now.

\ref{preimage3-i}: The structure of $\pF[J_2]$ is clear and directly
follows from Lemma~\ref{lem:preimage2}\ref{preimage2-b} (note $\pmSg{\id^{\pm_{1}}}=\pmSg{\id^{+}}\subsetneqq\pmSg{\id^{\pm_{2}}}$).

\ref{preimage3-ii}: Each $\pm$\hyp{}preclone $P$ in $\pF_{0}$ must contain some
$c_{0}^{\lambda}$ (for some signum $\lambda$) and for $f\in P$ we
must have $\mathring f\in F_{0}$ (including $\mathring f=\id$). Taking
into account Lemma~\ref{lem:preimage2}, all joins of the
join\hyp{}irreducibles give the sublattice $\pF[J_2]\cup\pF_{0}$ as
presented in Figure~\ref{fig:preimagesA}. The only ``nontrivial'' part
which needs an additional proof is the $\pm$\hyp{}preclone $P_{F_{0}}$ with
three lower covers $\pmSg{c_{0}^{-}}$,
$\pmSg{c_{0}^{+},\id^{\pm_{2}}}$ and $P_{J_{2}}=\pmSg{\id^{-}}$. This
case can be reduced to  proving
$\pmSg{c_{0}^{\pm},\id^{-}}=\pmSg{c_{0}^{-},\id^{\pm_{2}}}$. In fact,
we have 
$\id^{\pm_{2}}\in\pmSg{\id^{-}}$ by \ref{preimage3-i}, 
$c_{0}^{-}(x)=c_{0}^{\pm}(\id^{-}(x),x)$, and
$\id^{-}(x)=\id^{\pm_{2}}(c^{-}_{0}(x),x)$,
$c^{\pm}_{0}(x,y)=c^{-}_{0}(c^{-}_{0}(\id^{\pm_{2}}(x,y)))$.

\ref{preimage3-iii}: This case is dual to \ref{preimage3-ii}
(interchange $0$ and $1$).

\ref{preimage3-iv}: It is straightforward to see that all joins of the
join\hyp{}irreducibles (generated by) $c^{\pm}_{0},c^{+}_{0}$, $c^{\pm}_{1},c^{+}_{1}$
produce four new $\pm$\hyp{}preclones and together with the join\hyp{}irreducible
$\id^{\pm_{2}}$ we get the eight $\pm$\hyp{}preclones visible as ``cube'' in
the part $\pF_{01}$ of
 Figure~\ref{fig:preimagesA}. The joins with the
remaining join\hyp{}irreducibles $c^{-}_{0},c^{-}_{1}$ and $\id^{-}$
produce only two new $\pm$\hyp{}preclones above the cube (this
easily follows from the special properties of these join\hyp{}irreducibles, in
particular that they are idempotent and can change signa via
composition).  

\ref{preimage3-v}: There are four join\hyp{}irreducibles in $\pF_{\lnot}$
and again it is straightforward to prove that they generate the semi-interval
$\pF_{\lnot}$ as given in Figure~\ref{fig:preimagesA} (mainly using
Lemma~\ref{lem:preimage2} for $g=\lnot$ and proving
$\pmSg{\lnot^{\pm_{1}}}\lor\pmSg{\lnot^{-}}=P_{F_{\lnot}}$). Since
$F_{\lnot}$ is minimal (e.g., it does not contain constants), the $\pm$\hyp{}preclones in
$\pF_{\lnot}$ can only contain $\pm$\hyp{}preclones from $\pF[J_{2}]$ and
it is easy to see that
$\id^{\pm_{2}}\in\pmSg{\lnot^{\pm_{1}}}\lor\pmSg{\lnot^{\pm_{2}}}$ 
(since $\id^{\pm_{2}}(x,y)=y=\lnot^{\pm_{2}}(x,\lnot^{\pm_{1}}(y,x))$)
while $\id^{-}$ is contained only in $P_{F_{\lnot}}$.

\ref{preimage3-vi}: The hardest (and most tedious) job is to prove the
correctness of the presentation of $\pF_{2}$ in
Figure~\ref{fig:preimagesA}. We checked many parts by hand but the
final result was done with a computer program. Here some hints: Each
$\pm$\hyp{}preclone in $\pF_{2}$ must be a join of join\hyp{}irreducibles
containing negation and both constants, 
therefore it must be a join $H\lor K$ with $H\in\pF_{\lnot}$ and
$K\in\pF_{01}$. It is easy to check that the join of
$\pmSg{\lnot^{\pm_{1}}}$ and the join of
$\pmSg{\lnot^{\pm_{2}}}$ with the ``cube'' (see \ref{preimage3-iv}
above) in $\pF_{01}$ gives two ``cubes'' in $\pF_{2}$ with the same
upper face, altogether twelve $\pm$\hyp{}preclones (``visible'' in the right
lower part of $\pF_{2}$ in Figure~\ref{fig:preimagesA}). The remaining seven
$\pm$\hyp{}preclones in $\pF_{2}$ are the six joins of 
$H_{1},H_{2}\in\pF_{\lnot}$ with 
$K_{1},K_{2},K_{3}\in\pF_{01}$ (explicitly marked in
Figure~\ref{fig:preimagesA}), and the join $\pmSg{\lnot^{-}}\lor
K_{2}$. 
\end{proof}

Now we consider the preimage classes of further (non\hyp{}unary)
minimal clones: let 
 $F_{\land}$, $F_{\lor}$, and $F_{\oplus}$  be the
 minimal clones generated 
by conjunction ($\land$), disjunction ($\lor$) and $g_{\oplus}$,  respectively. 
Here $g_{\oplus}$ is the ternary Boolean minority 
function defined by $g_{\oplus}(x,y,z):=x\oplus y\oplus z$
($=\mathring{f}_{\oplus}(x,y,z,0)$, cf.\ \ref{min-theorem}\ref{min-thmA})

The
preimage classes are given in the following proposition.

Recall that, for a Boolean function $f:A^{n}\to A$, we denote by
$f^{\lambda}$ the corresponding $\pm$\hyp{}operation with signum
$\lambda\in\{+,-\}^{n}$ (we have $(f^{\lambda})\mathring{\,}=f$). 
The operators $\nabla^{+}$ and $\nabla^{-}$ add a fictitious variable
of sign $+$ or $-$, respectively,
at the first place (cf.\ Definition~\refI{2.3(4)}), e.g.,
$\nabla^{-}(\land^{(+,+)})=(\nabla(\land))^{(-,+,+)}$.
Note that the $\pm$\hyp{}preclones $\pmSg{\nabla^{-}(\land^{(+,+)})}$ and 
$\pmSg{\nabla^{-}(g_{\oplus}^{(+,+,+)})}$ are minimal $\pm$\hyp{}preclones
since $\pmSg{\nabla^{-}(\land^{(+,+)})}=\pmSg{g_{0000}^{(+,+,-)}}$ (cf.\
Theorem~\ref{min-theorem}\ref{min-thmC}), and
$\pmSg{\nabla^{-}(g_{\oplus}^{(+,+,+)})}=\pmSg{f_{\oplus}}$ (cf.\
Theorem~\ref{min-theorem}\ref{min-thmA}).

\normalsize
\begin{proposition}\label{prop:preimage4}
\leavevmode
  \begin{enumerate}[label={\upshape(\alph*)}]
  \item\label{preimage4a} 
$\pF[F_{\land}]=
\{\pmSg{\nabla^{-}\land^{(+,+)}},
\pmSg{\land^{(+,-)}},
\pmSg{\land^{(+,+)}},
\pmSg{\nabla^{+}\land^{(-,-)}},
\pmSg{\land^{(-,-)}},
\\ 
\hspace*{7ex}
\pmSg{\nabla^{-}\land^{(+,+)},\land^{(+,-)}},
\pmSg{\land^{(+,+)},\land^{(+,-)}},
\pmSg{\land^{(+,+)},\id^{\pm_{2}}}\}$
\hspace*{\fill} $|\pF[F_{\land}]|=8$,

\item\label{preimage4b} 
$\pF[F_{\lor}]=
\{\pmSg{\nabla^{-}\lor^{(+,+)}},
\pmSg{\lor^{(+,-)}},
\pmSg{\lor^{(+,+)}},
\pmSg{\nabla^{+}\lor^{(-,-)}},
\pmSg{\lor^{(-,-)}},
\\ 
\hspace*{7ex}
\pmSg{\nabla^{-}\lor^{(+,+)},\lor^{(+,-)}},
\pmSg{\lor^{(+,+)},\lor^{(+,-)}},
\pmSg{\lor^{(+,+)},\id^{\pm_{2}}}\}$
\hspace*{\fill} $|\pF[F_{\lor}]|=8$,

\item\label{preimage4c}
$\pF[F_{\oplus}]=
\{
\pmSg{\nabla^{-}g_{\oplus}^{(+,+,+)}},
\pmSg{g_{\oplus}^{(+,-,-)}},
\pmSg{g_{\oplus}^{(+,+,-)}},
\pmSg{g_{\oplus}^{(+,+,+)}},
\pmSg{g_{\oplus}^{(-,-,-)}},\\\hspace*{7ex}
\pmSg{g_{\oplus}^{(+,+,+)},g_{\oplus}^{(+,-,-)}},
\pmSg{g_{\oplus}^{(+,+,+)},g_{\oplus}^{(+,+,-)}}
\}$
 \hspace*{\fill} $|\pF[F_{\oplus}]|=7$.

  \end{enumerate}
\end{proposition}

\begin{figure}
  \begin{center}
\includegraphics{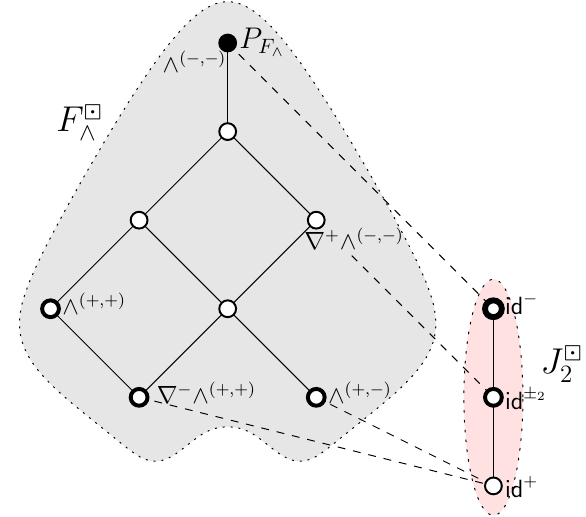}
  \end{center}
  \caption{The preimage class $\pF_{\land}$
     and its connection to $\pF[J_2]$\label{fig:preimagesB} 
(analogously for $\pF_{\lor}$)}
\end{figure}

\begin{proof}
\ref{preimage4a}: The five listed one-generated $\pm$\hyp{}preclones are different
and satisfy the following
inclusions:
$\pmSg{\nabla^{-}\land^{(+,+)}}\subset\pmSg{\land^{(+,+)}}
     \subset\pmSg{\land^{(-,-)}}$ and
$\pmSg{\land^{(+,-)}}\subset\pmSg{\nabla^{+}\land^{(-,-)}}
      \subset\pmSg{\land^{(-,-)}}$ (this is easy to see, e.g.,
      $\pm$\hyp{}operations with mixed signum can never generate a
      $\pm$\hyp{}operation with pure signum; and mixed or pure positive signum never
      can generate pure negative signum). It is straightforward to see
      that the joins of these $\pm$\hyp{}preclones give the ``lattice'' in
      Figure~\ref{fig:preimagesB}.

Further, an $n$\hyp{}ary Boolean
function $f\in F_{\land}$ ($n\geq 2$) is of the form
$f(x_{1},\dots,x_{n})=x_{1}\land \ldots\land x_{n}$ (up to fictitious variables).
It turns out that, for an arbitrary signum $\lambda\in\{+,-\}^{n}$, if $f$ has no fictitious variables, then the corresponding $\pm$\hyp{}operation $f^{\lambda}$
generates one of the listed (one-generated) $\pm$\hyp{}preclones, which
completes the proof.
 In fact we have:

$\pmSg{f^{\lambda}}=\pmSg{\land^{(+,+)}}$ for $\lambda=(+,\ldots,+)$
(pure positive signum),

$\pmSg{f^{\lambda}}=\pmSg{\land^{(+,-)}}$ for
$\lambda=(+,\ldots,+,-,\ldots,-)$ (mixed signum),

$\pmSg{f^{\lambda}}=\pmSg{\land^{(-,-)}}$ for $\lambda=(-,\ldots,-)$
(pure negative signum).

For taking into account fictitious variables it is enough to consider
$\pm$\hyp{}operations of pure signum (adding a fictitious variable to a
$\pm$\hyp{}operation of mixed signum can be undone by identification with another variable of the same sign).
Let $\lambda=(+,\dots,+,-,\dots,-)$ (with $m$ occurrences of `$+$', $1\leq m<n$).
It remains to check the following $\pm$\hyp{}operations
$f(x_{1},\dots,x_{m},\dots,x_{n}):=x_{1}\land\ldots\land x_{m}$ and 
$h(x_{1},\dots, x_{m+1},\dots,x_{n})=x_{m+1}\land\ldots\land x_{n}$.
But then we have
$\pmSg{f^{\lambda}}=\pmSg{\nabla^{-}\land^{(+,+)}}$ 
and $\pmSg{h^{\lambda}}=\pmSg{\nabla^{+}\land^{(-,-)}}$, which completes
the proof of \ref{preimage4a}.

\ref{preimage4b}: This can clearly be proved analogously to \ref{preimage4a}
by the duality of $\land$ and $\lor$.

\ref{preimage4c}: 
Each $n$\hyp{}ary term operation $h\in F_{\oplus}$ has an odd number of
essential variables (and it is at least ternary if
nontrivial). Identification of two variables makes them 
fictitious (since $x\oplus x=0$, $y+0=y$). Let
$h_{n}(x_{1},\dots,x_{n}):=x_{1}\oplus\ldots\oplus x_{n}$ and
$\lambda=(\underbrace{+,\ldots,+}_{m},-,\dots,-)$ with $0\leq m\leq n$, $3\leq n$.
Then we have:

$\pmSg{h_{n}^{\lambda}}=\pmSg{g_{\oplus}^{(+,-,-)}}$ if $m\neq n$ is odd,\\
$\pmSg{h_{n}^{\lambda}}=\pmSg{g_{\oplus}^{(+,+,-)}}$ if $m\neq 0$ is even,\\
$\pmSg{h_{n}^{\lambda}}=\pmSg{g_{\oplus}^{(+,+,+)}}$ if $m=n$,\\
$\pmSg{h_{n}^{\lambda}}=\pmSg{g_{\oplus}^{(-,-,-)}}$ if $m=0$.

This can be checked easily using identification of variables and
composition, e.g.,
$g_{\oplus}^{(+,-,-)}(x_{1},x_{2},x_{3})
             =h_{5}^{(+,+,+,-,-)}(x_{1},x_{1},x_{1},x_{2},x_{3})$,
$h_{5}^{(+,+,+,-,-)}(x_{1},\dots,x_{5})=g_{\oplus}^{(+,-,-)}(x_{1},x_{4},g_{\oplus}^{(+,-,-)}(x_{5},x_{2},x_{3}))$.

Thus each $n$\hyp{}ary term operation without fictitious variables
generates one of the above $\pm$\hyp{}preclones.
Taking fictitious variables into account, one has to add
$\pmSg{\nabla^{-}g_{\oplus}^{(+,+,+)}}$ and
$\pmSg{\nabla^{+}g_{\oplus}^{(-,-,-)}}$. However
$\pmSg{\nabla^{+}g_{\oplus}^{(-,-,-)}}=\pmSg{g_{\oplus}^{(+,+,-)}}$
since with the notation 
$h:=\nabla^{+}g_{\oplus}^{(-,-,-)}$ (note $\sgn(h)=(+,-,-,-)$) we have:
\begin{align*}
  h(x_{1},x_{2},x_{3},x_{4})
    &=g_{\oplus}^{(-,-,-)}(x_{2},x_{3},x_{4})\\
  &=g_{\oplus}^{(+,+,-)}(x_{1},g_{\oplus}^{(+,+,-)}(x_{1},x_{1},x_{2}),
                         g_{\oplus}^{(+,+,-)}(x_{3},x_{4},x_{1})),\\
 g_{\oplus}^{(+,+,-)}(x_{1},x_{2},x_{3})
 &=h(x_{1},h(x_{4},h(x_{1},x_{3},x_{3},x_{3}),x_{1},x_{2}),x_{3},x_{3}).
\end{align*}

Consequently each element of
the preimage class $\pF_{\oplus}$ is one of the listed one-generated
$\pm$\hyp{}preclones or the join of some of them. The final result can be seen
in Figure~\ref{fig:preimagesC}. The inclusions are strict and can be
proved again 
using composition;
to illustrate this, we present the proof for a single representative case, namely
$\pmSg{g_{\oplus}^{(+,-,-)}}\subset
\pmSg{g_{\oplus}^{(+,+,-)}}$. In fact,
\begin{align*}
  g_{\oplus}^{(+,-,-)}(x_{1},x_{2},x_{3})
  =g_{\oplus}^{(+,+,-)}(x_{1},x_{1},g_{\oplus}^{(+,+,-)}(x_{2},x_{3},x_{1}))
  \quad=x_{1}\oplus x_{2}\oplus x_{3}.
\end{align*}
Furthermore, each $\pm$\hyp{}operation in $\pmSg{g_{\oplus}^{(+,-,-)}}$
has an odd number of essential variables with sign `$+$' and
therefore can never generate $g_{\oplus}^{(+,+,-)}$.

As
Figure~\ref{fig:preimagesC} shows, there are only
two more elements in $\pF_{\oplus}$ which are the join of two
one-generated $\pm$\hyp{}preclones (namely of $g_{\oplus}^{(+,+,+)}$ with
$g_{\oplus}^{(+,-,-)}$ and $g_{\oplus}^{(+,+,-)}$, respectively). We
omit further details of the (straightforward and tedious) proof.

\begin{figure}
  \begin{center}
\includegraphics{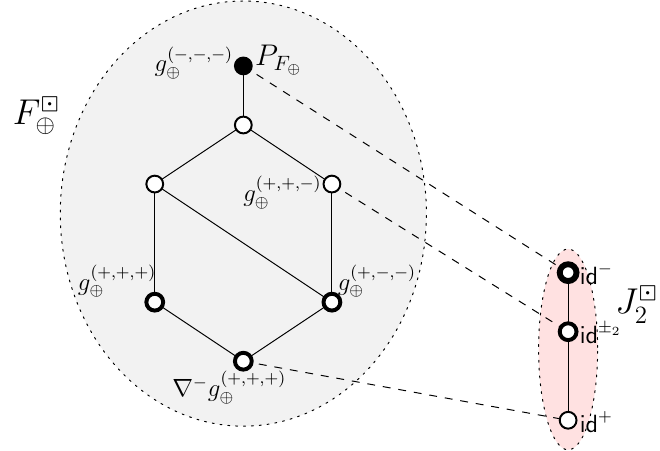}
  \end{center}
  \caption{The preimage class $\pF_{\oplus}$
     and its connection to $\pF[J_2]$\label{fig:preimagesC} }
\end{figure}
\end{proof}

There is one more minimal Boolean clone not considered yet, namely the
 minimal clone $F_{\mathrm{maj}}$ generated by the majority
 operation $m(x,y,z):=(x\land y)\lor (y\land z)\lor (z\land x)$. 
The structure of the preimage class $\pF_{\mathrm{maj}}$ is still not
fully described, but clearly $\pmSg{f_{\mathrm{maj}}}$ belongs to this
preimages class ($f_{\mathrm{maj}}$ equals $m$ but has an additional
fictitious variable of sign `$-$'; see
Theorem~\ref{min-theorem}\ref{min-thmA}). 

The Baker--Pixley theorem shows that
the polymorphism clone of a $(d+1)$\hyp{}ary near\hyp{}unanimity operation $g$ is
determined by its $d$\hyp{}ary invariant
relations and, consequently, that there are only finitely many clones
containing $g$. This can be generalized to $\pm$\hyp{}clones. We formulate
it here only for the special case of Boolean $\pm$\hyp{}clones containing
the majority $\pm$\hyp{}operation $f_{\mathrm{maj}}$ (thus $d=2$, a ternary
near\hyp{}unanimity operation is a majority operation). However, the proof remains
valid if $A$ is an arbitrary finite set and $f_{\mathrm{maj}}$ is an
arbitrary majority  operation with an additional fictitious variable 
with sign `$-$'. Moreover, it is straightforward to generalize this
result to $(d+1)$\hyp{}ary near\hyp{}unanimity operations with a single negative
sign for an additional fictitious variable.

\begin{proposition}[Baker--Pixley theorem for Boolean
  $\pm$\hyp{}preclones]\label{prop:BakerPixley}
Let $P$ be a $\pm$\hyp{}preclone containing $f_{\mathrm{maj}}$.
Then $P$ is determined by its binary $\pm$\hyp{}relations:
$P=\pmPol\pmInva[2]P$. 
In particular, the interval
$[\pmSg{f_{\mathrm{maj}}},\pmOp(A)]_{\pmL_{2}}$ 
is finite.
\end{proposition}

\begin{proof} Clearly $P\subseteq\pmPol\pmInva[2]P$. Thus only the
  inverse containment must be proved. Let
  $f^{\lambda}\in \pmPola[\lambda]\pmInva[2]P$ with signum
  $\lambda=(s_{1},\dots,s_{n})\in \{+,-\}^{n}$. We are going to show
  $f^{\lambda}\in P$, which will complete the proof. Assume that $\lambda$
  contains at least one sign $-$ (the case of pure positive signum
  will be discussed later), say $s_{i_{0}}=-$.
Let 
\begin{align*}
  \cM:=\{B\subseteq A^{n}\mid\exists g_{B}\in
  P^{(\lambda)}:
f\restriction_{B}=g_{B}\restriction_{B}\}.
\end{align*}
We show $A^{n}\in\cM$ by induction on the cardinality $t:=|B|$. Then
$f=g_{A^{n}}\in P$ and we are done. 

Let $t=2$, i.e., $B$ is of the form
$B=\{(a_{1},\dots,a_{n}),(b_{1},\dots,b_{n})\}$.

Consider the invariant
(for $P$) $\pm$\hyp{}relation $\Gamma_{P}(\rho)$ (cf.\ Definition~\refI{3.11}) generated by 
$\rho:=\left(
  \begin{smallmatrix}
    a_{1}&\ldots&a_{n}\\b_{1}&\ldots&b_{n}
  \end{smallmatrix}\right)$ with
$\blambda_{\rho}:=\lambda$ (for the notation, see Definition~\refI{2.8}).
Since $\Gamma_{P}(\rho)\in\pmInva[2]P$ we have
$f\pmpreserves\Gamma_{P}(\rho)$, in particular
$(a,b):=(f(a_{1},\dots,a_{n}),\allowbreak f(b_{1},\dots,b_{n}))\in
(\Gamma_{P}(\rho))_{+}$. 

From Proposition~\refI{3.13(iii)} we get
$(\Gamma_{P}(\rho))_{+}=\{g(r_{1},\dots,r_{n})\mid g\in 
    P^{(\blambda_{\rho})}\}$ where 
$r_{i}=\left(\begin{smallmatrix}  a_{i}\\b_{i} \end{smallmatrix}\right)$
(note that here $+$ is the unit element $e$ of $S=\pm$).
Therefore there exists 
a $\pm$\hyp{}operation $g_{B}\in P^{(\lambda)}$ (note
$\lambda=\blambda_{\rho}$ by construction) such that 
$(a,b)=(g_{B}(a_{1},\dots,a_{n}),g_{B}(b_{1},\dots,b_{n}))$, 
i.e., $f\restriction_{B}=g_{B}\restriction_{B}$.

(Note that the above argument works also for $t = 1$; in this case we have $(a_1, \dots, a_n) = (b_1, \dots, b_n)$.)

Now assume that $B'\in\cM$ for all $B'$ with $t$ elements 
($2\leq t\leq |A^{n}|-1$). We show that $B\in\cM$ for all $B$ with $t+1$
elements. Then, by induction on $t$, $A^{n}\in\cM$ and we are done.

Let $B\subseteq A^{n}$ with $|B|=t+1$. Choose three different elements
$\bb_{1},\bb_{2},\bb_{3}\in B$ (this is possible because $|B|=t+1\geq
3$) and let $B_{i}:=B\setminus\{\bb_{i}\}$ ($i=1,2,3$). Then
$|B_{i}|=t$ and therefore $B_{i}\in\cM$ by assumption. Thus there
exist $g_{B_{i}}\in P^{(\lambda)}$ 
satisfying $f(\bb)=g_{B_{i}}(\bb)$
for all $\bb\in B_{i}$ ($i=1,2,3$).

Then the $\pm$\hyp{}operation $g$ of signum $\lambda$ (recall that $s_{i_0} = \mathord{-}$ by assumption) given by the rule
\[
g(\bx) = f_{\mathrm{maj}}(g_{B_{1}}(\bx), g_{B_{2}}(\bx), g_{B_{3}}(\bx), x_{i_{0}})
\quad
\text{for all $\bx \in A^n$}
\]
is in $P$, and we have $f(\bb) = g(\bb)$ for all $\bb \in B$.
In fact, for $\bb = \bb_i$ with $\{i, i', i''\} = \{1, 2, 3\}$, we have $g_{B_{i'}}(\bb) = g_{B_{i''}}(\bb) = f(\bb)$ (since $\bb_{i} \in B_{j}$ for $j \neq i$), and hence $g(\bb) = f(\bb)$ because $f_{\mathrm{maj}}$ is a majority operation on the first three arguments and the fourth argument is fictitious.
For $\bb \in B \setminus \{\bb_1, \bb_2, \bb_3\}$, we trivially have $g_{B_1}(\bb) = g_{B_2}(\bb) = g_{B_3}(\bb) = f(\bb)$, and again $g(\bb) = f(\bb)$.
This shows $f\restriction_{B} = g \restriction_{B}$.

Finally we have to consider the case $\sgn(f)=(+,\dots,+)$.
If $f$ is a projection, then $f \in P$ and we are done, so assume that $f$ is not a projection.
Then $\mathring f\in F_{\mathrm{maj}}=\Sg{m}$ and therefore (because
$F_{\mathrm{maj}}$ is a minimal clone) we have $m\in\Sg{\mathring f}$, which
implies $m^{+++}\in\pmSg{f}\subseteq P$ (note that the positive part
$P^{[+]}$ of $P$ is a clone, cf.\ Proposition~\refI{4.7}, here 
$\mathring P^{[+]}=F_{\mathrm{maj}}$). Now we can proceed as above (as for the
classical Baker--Pixley theorem) using $m^{+++}$ instead of
$f_{\mathrm{maj}}$:
\[
f(\bb)=m^{+++}(g_{B_{1}}(\bb),g_{B_{2}}(\bb),g_{B_{3}}(\bb)),
\]
which (by induction) implies $f\in P$, and we are done.
\end{proof}


\section{Concluding remarks}

We described all 9 coatoms and 20 atoms of the lattice $\pmL_2$ of Boolean $\pm$-preclones, and obtained an analogous result for the sublattice $\pmLsd_2$ of self-dual Boolean $\pm$-preclones. In contrast to Post's lattice of Boolean clones, $\pmL_2$ is uncountable. Boolean $\pm$-preclones are grouped into preimage classes of a corresponding Boolean clone, and we described the preimage classes of all unary Boolean clones as well as of the minimal Boolean clones, except for the majority clone. A Baker--Pixley type theorem for $\pm$-preclones containing the majority operation shows that certain subintervals of $\pmL_2$ remain finite and well-behaved. These results provide some structural insight into the lattice $\pmL_2$, and we mention here some further questions.

\begin{itemize}
\item Which of the ideals in $\pmL_2$ below maximal Boolean $\pm$-preclones
are countable, and which are uncountable?
\item Find the cardinality of the sublattice $\pmLsd_2$ of self-dual $\pm$-preclones.
While we have classified its maximal and minimal elements, the overall structure is yet to be
determined.
\item Describe the preimage class of the majority clone $F_{\mathrm{maj}}$.
\item While we have focused on the Boolean case $A = \{0,1\}$, some of our results, such
as the maximality of $\pmPol(\leq,\geq)$ for finite posets and the
Baker--Pixley theorem for majority operations, hold for arbitrary
finite base sets.
Which other results remain valid for arbitrary base sets?
\item Extend the classification of maximal $\pm$-preclones to general finite sets $A$.
\end{itemize}

Finally, we refer to the list of open problems mentioned in Section~\refI{6} in Part~I.

\def\cprime{$'$} \def\cprime{$'$}

\end{document}